%% file: ncube4fin.tex
\def\version{December 22, 2003}

\documentclass[12pt]{article}
\usepackage[psamsfonts]{amsfonts}
\usepackage{amsmath,amssymb,amsthm}
\usepackage[dvips]{graphicx}
\usepackage{eepic}
\input def99c % for LaTeX2e, enables \eq  and \en (!!)
\UseSection  %Necessary to define Numbering scheme for Theorem, etc.
\renewcommand{\to}      {\rightarrow}

\newcommand{\Pro}{{\mathbb P}_p}

\newcommand{\ben}{\begin{enumerate}}
\newcommand{\een}{\end{enumerate}}

\newcommand{\xleftrightarrow}[2][]{\leftarrow\hspace{-2.0ex}\xrightarrow[#1]{#2}}

\newcommand{\xconn}[1]{\xleftrightarrow[#1]{}}

\newcommand{\shift}   {\!\!\!\!}
\newcommand{\SSS}   {\sss}

\newcounter{countC}  % Defined the counter ``countC''
\newcounter{countR}  % Defined the counter ``countR''
\newcommand{\Z}{\Zbold}

\newcommand{\conn}{\leftrightarrow}

\newcommand{\dbc}{\Leftrightarrow}
\newcommand{\ct}[1]     { \stackrel{#1}{\conn} }

\newcommand{\smallsup}[1] {{\scriptscriptstyle{({#1}})}}

\newcommand{\bigo}{O}

\newcommand{\gr}{\mathbb G}
\newcommand{\qn}{{\mathbb Q}_n}

\newcommand{\cn}{\Omega}

\title  {
        Expansion in $n^{-1}$ for
        percolation critical values
        \\
        on the $n$-cube and $\Z^n$:
        the first three terms
        }

\author{
Remco van der Hofstad\thanks{Department of Mathematics and Computer Science,
Eindhoven University of Technology, P.O.\ Box  513,
5600 MB Eindhoven, The Netherlands.
{\tt rhofstad@win.tue.nl}}
\and
Gordon Slade\thanks{Department of Mathematics, University of British Columbia,
Vancouver, BC V6T 1Z2, Canada. {\tt slade@math.ubc.ca}}
}

\date\version

\begin{document}

\maketitle

\begin{abstract}
Let $p_c(\qn)$ and $p_c(\Z^n)$ denote the critical values for
nearest-neighbour bond percolation on
the $n$-cube $\qn = \{0,1\}^n$ and on $\Z^n$, respectively.
Let $\cn = n$ for $\gr = \qn$ and $\cn = 2n$ for $\gr = \Z^n$
denote the degree of $\gr$.
We use the lace expansion to prove that for both $\gr = \qn$
and $\gr = \Z^n$,
\eqalign
    p_c(\gr)  & =
    \cn^{-1} + \cn^{-2} + \frac{7}{2} \cn^{-3}
    + O(\cn^{-4}).
    \nonumber
\enalign
This extends by two terms the
result $p_c(\qn) = \cn^{-1} + O(\cn^{-2})$
of Borgs, Chayes, van der Hofstad, Slade and Spencer,
and provides a simplified proof of
a previous result of Hara and Slade for $\Z^n$.
\end{abstract}

%%%%%%%%%%%%%%%%%%%%%%%%%%%%%%%%%%%%%%%%%%%%%%%%%%%%%%%%%%%%%%%%%%%%%%%%%%%%%%%
%%%%%%%%%%%%%%%%%%%%%%%%%%%%%%%%%%%%%%%%%%%%%%%%%%%%%%%%%%%%%%%%%%%%%%%%%%%%%%%
\section{Main result}
\label{sec-intro}

We consider bond percolation on $\Z^n$ with edge set consisting
of pairs $\{x,y\}$ of vertices in $\Z^n$ with $\|x-y\|_1 = 1$,
where $\|w\|_1 = \sum_{j=1}^n |w_j|$ for $w \in \Z^n$.
Bonds (edges) are independently occupied with probability $p$
and vacant with probability $1-p$.
We also consider bond percolation on the $n$-cube $\qn$, which
has vertex set $\{0,1\}^n$ and edge set consisting
of pairs $\{x,y\}$ of vertices in $\{0,1\}^n$
with $\|x-y\|_1 = 1$, where we regard
$\qn$ as an additive group with addition component-wise
modulo~2.
Again
bonds are independently occupied with probability $p$
and vacant with probability $1-p$.
We write $\gr$ in place of $\qn$ and $\Z^n$ when we wish
to refer to both models simultaneously.
We write $\cn$ for the degree of $\gr$, so that
$\cn = 2n$ for $\Z^n$ and $\cn = n$ for $\qn$.

For the case of $\Z^n$, the critical value is defined by
\eq
    p_c(\Z^n) = \inf \{p : \exists \text{  an infinite
    connected cluster of occupied bonds a.s.}\}.
\en
Given a vertex $x$ of $\gr$, let $C(x)$ denote the connected
cluster of $x$, i.e., the set of vertices $y$ such that
$y$ is connected
to $x$ by a path consisting of occupied bonds.  Let $|C(x)|$
denote the cardinality of $C(x)$,
and let $\chi(p) = \Ebold_p|C(0)|$ denote the expected cluster size of
the origin.
Results of \cite{AB87,Mens86} imply that
\eq
    p_c(\Z^n) = \sup\{p : \chi(p) <\infty\}.
\en
is an equivalent definition of the critical value.

For percolation on a finite graph $\gr$, such as $\qn$, the above
characterizations of $p_c(\gr)$ are inapplicable. In
\cite{BCHSS04a,BCHSS04b,BCHSS04c} (in particular, see
\cite{BCHSS04c}), it was shown that there is a small positive
constant $\lambda_0$ such that the critical value $p_c(\qn) =
p_c(\qn;\lambda_0)$ for the $n$-cube is defined implicitly by
\eq
\lbeq{pcndef}
    \chi(p_c(\qn)) = \lambda_0 2^{n/3}.
\en
Given $\lambda_0$, \refeq{pcndef} uniquely specifies
$p_c(\qn)$, since $\chi(p)$ is a polynomial in $p$ that increases
from $\chi(0)=1$ to $\chi(1)=2^n$.
%As discussed in
%\cite{BCHSS04c}, $p_c(\qn;\lambda_0)$ depends only weakly on the
%choice of $\lambda_0$.  This point is also reflected in the Remark
%below.

Our main result is the following theorem.

\begin{theorem}
\label{thm-mr}
{\rm (i)}
For $\gr=\Z^n$,
\eq
\lbeq{mrZn}
    p_c(\Z^n) = \frac{1}{2n} + \frac{1}{(2n)^2}
    + \frac{7}{2} \frac{1}{(2n)^3} + O\big(\frac{1}{(2n)^4}\big)
    \quad
    \quad
    \mbox{as $n \to \infty$}.
\en
{\rm (ii)}
For $\qn$, fix constants $c,c'$
independent of $n$, and choose $p$ such that
$\chi(p) \in [cn^{3}, c'n^{-6}2^{n}]$ (e.g., $p=p_c(\qn;\lambda_0)$).
Then
\eq
\lbeq{mrqn}
    p = \frac{1}{n} + \frac{1}{n^2}
    + \frac{7}{2} \frac{1}{n^3} + O\big(\frac{1}{n^4}\big)
    \quad
    \quad
    \mbox{as $n \to \infty$}.
\en
The constant in the error term depends on $c,c'$, but does not
depend otherwise on $p$.
\end{theorem}

By Theorem~\ref{thm-mr},
the expansions of $p_c(\gr)$ in powers of $\cn^{-1}$
are the same for $\qn$ and $\Z^n$, up to
and including order $\cn^{-3}$.
Higher order coefficients could be computed using our methods, but
the labour cost increases sharply with each subsequent term.
Although we stop short of computing the coefficient
of $\cn^{-4}$, we
%mention indications below that make us
expect that the coefficients
for $\qn$ and $\Z^n$ will differ at this order.
In \cite{HS04b},
for both $\qn$ and $\Z^n$, we prove the existence of
asymptotic expansions for $p_c(\gr)$ to all orders in $\cn^{-1}$,
without computing the numerical values of the coefficients.

For $\qn$, it was shown by Ajtai, Koml\'os and Szemer\'edi \cite{AKS82}
that $p_c(\qn) > n^{-1}(1+\epsilon)$ for every {\em fixed} $\epsilon >0$
(although the above definition of $p_c(\qn)$ did not appear until
\cite{BCHSS04a}).
Bollob\'as, Kohayakawa and {\L}uczak \cite{BKL92}
improved this to $p_c(\qn) \in [\frac{1-e^{-o(n)}}{n-1}, \frac{1}{n} +
60\frac{(\log n)^3}{n^2}]$.
Theorem~\ref{thm-mr} extends the very recent result
$p_c(\qn) = n^{-1} + O(n^{-2})$ of
\cite{BCHSS04a,BCHSS04b} by two terms.
Bollob\'as, Kohayakawa and {\L}uczak \cite{BKL92}
raised the question of whether
the critical value might be \emph{equal} to $\frac{1}{n-1}$, but we
see from \refeq{mrqn} that $p_c(\qn) = \frac{1}{n-1} +
\frac{5}{2} n^{-3} + O(n^{-4})$.

For $\Z^n$, Theorem~\ref{thm-mr} is
identical to a result
of Hara and Slade \cite{HS93up,HS95}.  Earlier,
Bollob\'as and Kohayakawa \cite{BK94},
Gordon \cite{Gord91}, Kesten \cite{Kest90} and Hara and Slade \cite{HS90a}
obtained the first term in \refeq{mrZn} for $\Z^n$ with error
terms $\bigo((\log n)^2n^{-2} )$,
$\bigo(n^{-65/64})$, $\bigo((\log \log n)^2(n\log n)^{-1})$
and $\bigo(n^{-2})$, respectively.  Recently, Alon, Benjamini and Stacey
\cite{ABS02}
gave an alternate proof that $p_c(\Z^n)$ is asymptotic to $(2n)^{-1}$
as $n \to \infty$.
The expansion
\eq
\lbeq{pcGR}
    p_c(\Z^n)
    =  \frac{1}{2n} + \frac{1}{(2n)^2} + \frac{7}{2(2n)^3}
    +\frac{16}{(2n)^4} + \frac{103}{(2n)^5} + \cdots
\en
was reported in \cite{GR78}, but with no rigorous bound on the remainder.

We remark that
for \emph{oriented} percolation on $\Z^n$, defined in such a way
that the forward degree is $n$, it was proved in
\cite{CD83} that the critical value obeys the bounds
\eq
    \frac{1}{n} + \frac{1}{2}\frac{1}{n^3} + o\Big( \frac{1}{n^3}\Big)
    \leq
    p_c(\text{oriented $\Z^n$})
    \leq
    \frac{1}{n} + \frac{1}{n^3} + O\Big( \frac{1}{n^4}\Big).
\en

Our method is based on the lace expansion and applies
the general approach of \cite{HS93up,HS95} that was used to prove
Theorem~\ref{thm-mr}(i) for $\Z^n$, but our method here is
simpler and applies to $\Z^n$ and $\qn$ simultaneously.

\noindent {\bf Remark.}
For $\qn$, it is a direct consequence of
\cite[Proposition~1.2]{HS04b} that
if there is some sequence $p$ (depending on $n$) with
$\chi(p) \in [cn^3,c'n^{-6}2^n]$ such that
$p = n^{-1} + n^{-2} + \frac{7}{2} n^{-3} + O(n^{-4})$,
then the same asymptotic formula holds for \emph{all}
such $p$.  Thus it suffices to prove \refeq{mrqn} for a single such
sequence $p$.
We fix some sequence $f_n$ such that $\lim_{n\to \infty} f_n n^{-M} = \infty$
for every positive integer $M$ and such that
$\lim_{n\to \infty} f_n e^{-\alpha n} = 0$ for every $\alpha >0$.
We define $\bar p$ by $\chi(\bar p) = f_n$, and observe that eventually
$\chi(\bar p) \in [cn^3,c'n^{-6}2^n]$.
For $\gr = \qn$, it therefore suffices to prove that
$\bar p$ has the expansion \refeq{mrqn}.
We will use the notation
\eq
\lbeq{pbardef}
    \bar p_c = \bar p_c(\gr) =
    \begin{cases}
    \bar p & (\gr = \qn),
    \\
    p_c(\Z^n) & (\gr = \Z^n).
    \end{cases}
\en

%%%%%%%%%%%%%%%%%%%%%%%%%%%%%%%%%%%%%%%%%%%%%%%%%%%%%%%%%%%%%%%%%%%%%%%%%%%%%%%
%%%%%%%%%%%%%%%%%%%%%%%%%%%%%%%%%%%%%%%%%%%%%%%%%%%%%%%%%%%%%%%%%%%%%%%%%%%%%%%
\section{Application of the lace expansion}
\label{sec-le}

For $\qn$ or $\Z^n$ with $n$ large, the lace expansion \cite{HS90a}
gives rise to an identity
        \eq
        \lbeq{tauk'}
        \chi(p) = \frac{1+\hat \Pi_p}
        {1-\cn p[1+\hat \Pi_p]},
        \en
where $\hat\Pi_p$ is a function
that is finite for $p \leq p_c(\gr)$.  Although we do not display
the dependence explicitly in the notation, $\hat\Pi_p$ does
depend on the graph $\qn$ or $\Z^n$.
The identity \refeq{tauk'} is valid for $p \leq p_c(\gr)$.
For a derivation of the lace expansion, see, e.g.,
\cite[Section~3]{BCHSS04b}.
It follows from \refeq{tauk'} that
\eq
\lbeq{pcchi}
    \cn p = \frac{1}{1+\hat \Pi_{p}}
    - \chi(p)^{-1} .
\en

The function $\hat\Pi_p$ has the form
    \eq
    \lbeq{2pt.37}
    \hat\Pi_p = \sum_{N=0}^{\infty} (-1)^N
    \hat\Pi_p^{\smallsup{N}},
    \en
with (recall \refeq{pbardef})
    \eq
    \lbeq{Pibds}
    |\hat \Pi_p^{\smallsup{N}}|
    \leq \left(\frac{C}{\cn }\right)^{N\vee 1}
    \quad
    \text{uniformly in $p \leq \bar p_c$}.
    \en
For $\qn$, the formula \refeq{tauk'} and the bounds \refeq{Pibds}
are given in \cite[(6.1)]{BCHSS04b} and \cite[Lemma~5.4]{BCHSS04b},
respectively (with our $\hat{\Pi}_p$ written as $\hat\Pi_p(0)$).
In more detail, \cite[Lemma~5.4]{BCHSS04b}
states that $\hat{\Pi}_p^\smallsup{N} \leq
[\mbox{const}(\lambda^3 \vee \beta)]^{N\vee 1}$, where
$\lambda = \chi(p) 2^{-n/3} \leq f_n 2^{-n/3}$ for
$p \leq \bar p_c(\qn)$.  By definition, $f_n2^{-n/3}$ is
exponentially small in $n$.
In addition,
it is shown in \cite[Proposition~2.1]{BCHSS04b} that
$\beta$ can be chosen proportional to $n^{-1}$.
It follows
from \refeq{pcchi} that
    \eq
    \lbeq{pcform}
    n \bar p_c(\qn)
    = \frac{1}{1+\hat \Pi_{\bar p_c(\qn)}}
    +O(f_n^{-1}) .
    \en
The second term on the right hand side of \refeq{pcform}
can be neglected in the proof of Theorem~\ref{thm-mr}.
Equations~\refeq{2pt.37}--\refeq{pcform}
give $\bar p_c(\qn) = n^{-1}+O(n^{-2})$.

For $\Z^n$, \refeq{tauk'} and \refeq{Pibds} follow from
results in \cite[Section~4.3.2]{HS90a}.
(Note the notational difference that
in \cite{HS90a} what we are calling here
$\hat{\Pi}_p^\smallsup{N}$ is called $\hat{g}_N(0)$
and that $\hat{\Pi}_p^\smallsup{N}$ in \cite{HS90a}
is something different.)
Since $\chi(p_c(\Z^n))=\infty$,
it follows
from \refeq{pcchi} that
\eq
\lbeq{pcdefZd}
    2np_c(\Z^n) =  \frac{1}{1+\hat \Pi_{p_c(\Z^n)}}.
\en
With \refeq{2pt.37}--\refeq{Pibds},
this implies that $p_c(\Z^n) = (2n)^{-1} + O(n^{-2})$.

The identities \refeq{pcform} and \refeq{pcdefZd} give recursive
equations for $\bar p_c$.
To prove Theorem~\ref{thm-mr} using this recursion, we will
apply the following proposition.
In its statement, we write
\eq
\lbeq{Omega'def}
    \cn' = \begin{cases} n - 1 & \text{for $\qn$} \\
        2n-2 & \text{for $\Z^n$}. \end{cases}
\en

%%%%%%%%%%%% PROP %%%%%%%%%%%%%%%%%%%%
    \begin{prop}
    \label{prop-Piasya}
    For $\gr = \Z^n$ and $\gr = \qn$, uniformly in $p\leq
    \bar p_c(\gr)$,
    \eqalign
    \lbeq{Pi0asy}
    \hat \Pi^\smallsup{0}_p
    &=\frac 32 \cn\cn' p^4 + \bigo(\cn^{-3}),\\
    \lbeq{Pi1asy}
    \hat \Pi^\smallsup{1}_p
    &=\cn p^2+4\cn\cn'p^4 +\bigo(\cn^{-3}),\\
    \lbeq{Pi2asy}
    \hat \Pi^\smallsup{2}_p
    &=\cn p^3+\cn (\cn-1) p^4  +\bigo(\cn^{-3}),\\
    \sum_{N=3}^\infty \hat \Pi^\smallsup{N}_p&=\bigo(\cn^{-3}).
    \lbeq{remainder}
    \enalign
    \end{prop}
%%%%%%%%%%%% PROP %%%%%%%%%%%%%%%%%%%%

We show now that Proposition~\ref{prop-Piasya} implies Theorem~\ref{thm-mr}.
It follows from
$\cn \bar p_c(\gr)= 1 + O(\cn^{-1})$
(as noted below \refeq{pcform} and \refeq{pcdefZd}),
\refeq{2pt.37}, and Proposition~\ref{prop-Piasya} that
\eq
\lbeq{Pilead}
    \hat\Pi_{\bar p_c(\gr)} = -\frac{1}{\cn} + \bigo(\cn^{-2}).
\en
With \refeq{pcform}--\refeq{pcdefZd}, this implies that
    \eq
    \cn \bar p_c(\gr) =1+\frac{1}{\cn}+\bigo(\cn^{-2}).
    \en
Using this
in the bounds of Proposition~\ref{prop-Piasya}, along
with  \refeq{2pt.37}, gives
    \eqalign
    \hat \Pi_{\bar p_c(\gr)} & =
    \frac{3}{2\cn^2} -\cn(\frac{1}{\cn}+\frac{1}{\cn^2})^2
    -\frac{4}{\cn^2} +\frac{1}{\cn^2} +\frac{1}{\cn^2}
    +\bigo(\cn^{-3})
    \nnb &
    =-\frac{1}{\cn}-\frac{5}{2\cn^2}+\bigo(\cn^{-3}).
    \enalign
Substitution of this improvement of \refeq{Pilead}
into \refeq{pcform}--\refeq{pcdefZd} then gives
    \eq
    \cn \bar p_c(\gr)
    %=\frac{1}{1-\frac1\cn-\frac{5}{2\cn^2}+\bigo(\cn^{-3})}
    =1+\frac{1}{\cn}
    +\frac{7}{2\cn^2}+\bigo(\cn^{-3}).
    \en
Thus, to prove Theorem~\ref{thm-mr}, it suffices to
prove Proposition~\ref{prop-Piasya}.  Since \refeq{remainder}
is a consequence of \refeq{Pibds}, we must prove
\refeq{Pi0asy}--\refeq{Pi2asy}.  Precise definitions of
$\hat\Pi^\smallsup{N}_p$, for $N=0,1,2$, will be given in
Section~\ref{sec-Pibds}.

%We note evidence in Proposition~\ref{prop-Piasya} that
%there can be a difference between the coefficients
%of $p_c(\Z^n)$ and $p_c(\qn)$ at
%order $\cn^{-4}$, arising from terms of order $\cn^{-3}$
%in the expansion for $\hat\Pi_p$.  For $\Z^n$, $\cn' = 2n-2$,
%while
%for $\qn$, $\cn'= n-1$.
%The factors $\cn'$ in \refeq{Pi0asy}--\refeq{Pi1asy} have
%combinatorial meaning, and the terms
%$-2$ and $-1$ in $\cn'$ in the terms
%$\frac 32 \cn\cn' p^4$ and $4\cn\cn'p^4$
%will have different
%effects at order $\cn^{-3}$.
%Although such differences
%will appear in both $\hat \Pi^\smallsup{0}_p$
%and $\hat \Pi^\smallsup{1}_p$, it can easily be checked
%that they do not cancel
%each other, and we do not expect exact cancellation
%from another source.

%%%%%%%%%%%%%%%%%%%%%%%%%%%%%%%%%%%%%%%%%%%%%%%%%%%%%%%%%%%%%
%%%%%%%%%%%%%%%%%%%%%%%%%%%%%%%%%%%%%%%%%%%%%%%%%%%%%%%%%%%%%
\section{Preliminaries}
\label{sec-pre}

Before proving Proposition~\ref{prop-Piasya}, we recall
and extend some estimates from \cite{BCHSS04b,HS90a}.

Let $D(x) = \cn^{-1}$ if $x$ is adjacent to $0$,
and $D(x)=0$ otherwise.
Thus $D(y-x)$ is the transition
probability for simple random walk on $\gr$ to make a step
from $x$ to $y$.  Let
$\tau_p(y-x) = \Pbold_p(x \conn y)$ denote the two-point
function.
For $i \geq 0$, we denote by
    \eq
    \lbeq{xconnidef}
    \{x\xconn{i} y\}
    \en
the event that $x$ is connected to $y$ by
an occupied (self-avoiding) path of length
at least $i$,
and define
\eq
    \tau_p^\smallsup{i}(x,y) = \Pbold(x\xconn{i} y).
\en

We define the Fourier transform of
an absolutely summable function $f$ on the vertex set $\Vbold$
of $\gr$ by
\eq
    \hat{f}(k) = \sum_{x \in \Vbold} f(x) e^{ik\cdot x}
    \quad \quad
    (k \in \Vbold^*),
\en
where $\Vbold^* = \{0,\pi\}^n$ for $\qn$ and $\Vbold^* = [-\pi,\pi]^n$
for $\Z^n$.
We write the inverse Fourier transform as
\eq
\lbeq{FTinv}
    f(x) = \int \hat{f}(k) e^{-ik\cdot x},
\en
where we use the convenient notation
\eq
    \int \hat{g}(k) = \begin{cases}
    2^{-n} \sum_{k \in \{0,\pi\}^n} \hat{g}(k) & (\gr = \qn)
    \\
    \int_{[-\pi,\pi]^n} \hat{g}(k) \frac{d^nk}{(2\pi)^n}
    & (\gr = \Z^n).
    \end{cases}
\en
Let
\eq
    (f*g)(x) = \sum_{y \in \Vbold} f(y)g(x-y)
\en
denote convolution, and
let $f^{*i}$ denote the convolution of $i$ factors of $f$.

Recall from \cite{AN84} that $\hat\tau_p(k) \geq 0$ for all $k$.
For $i,j$ non-negative integers, let
\eqalign
    T_p^\smallsup{i,j}
    & =
     \int  |\hat{D}(k)|^i \hat{\tau}_p(k)^j,
    \\
    \lbeq{Tp-def}
    T_p &= \sup_{x} (p\cn)(D*\tau_p^{*3})(x).
\enalign
We will use the following lemma, which provides minor
extensions of results of \cite{BCHSS04b,HS90a}.  The lemma
will also be useful in \cite{HS04b}.

\begin{lemma}
\label{lem-Tbd}
For $\gr = \Z^n$ and $\gr = \qn$,
there are constants $K_{i,j}$ and $K$
such that
for all $p \leq \bar p_c(\gr)$,
\eqalign
\lbeq{Tpij}
     T_p^\smallsup{i,j}
     & \leq K_{i,j} \cn^{-i/2} \quad \text{($i,j\geq 0$)},
     \\
\lbeq{Tp}
     T_p & \leq K \cn^{-1},
\\
\lbeq{tauiT}
    \sup_x \tau_p^\smallsup{i}(x)
    & \leq
    \begin{cases}
    K \cn^{-1} & (i =1) \\
    2^i K_{i,1}\cn^{-i/2} & (i \geq 2).
    \end{cases}
\enalign
The above bounds are valid for $n \geq 1$ for $\qn$, and for $n$ larger
than an absolute constant for $\Z^n$, except \refeq{Tpij} also requires
$n \geq 2j+1$ for $\Z^n$.
\end{lemma}

\proof
We prove the bounds \refeq{Tpij}--\refeq{tauiT} in sequence.

\smallskip \noindent {\em Proof of \refeq{Tpij}.}
We first prove that for $\Z^n$ and $\qn$, and for positive integers $i$,
there is a positive $a_i$ such that
\eq
\lbeq{retpr}
    \int \hat{D}(k)^{2i}
    \leq \frac{a_i}{\cn^i}.
\en
The left side is equal to the probability that a
random walk that starts at the origin returns to
the origin after $2i$ steps, and therefore is equal to
$\Omega^{-2i}$ times the number of walks that make the transition
from 0 to 0 in $2i$ steps. Each such walk must take an even number
of steps in each coordinate direction, so it must lie within
a subspace of dimension $\ell\leq \min\{i,n\}$.  If we fix the
subspace, then each step
in the subspace can be chosen from at most $2\ell$ different directions
(for $\qn$, from $\ell$ directions).
Thus, there are at most $(2\ell)^{2i}$ walks in the
subspace.  Since the number of subspaces of fixed
dimension $\ell$ is given by ${n\choose\ell}\leq n^\ell/\ell!$, we
obtain the bound
\eq
    \sum_{\ell=1}^{i}\frac 1{\ell!}n^\ell (2\ell)^{2i}
    \leq
    n^i i^{2i} \sum_{\ell=1}^{i}\frac 1{\ell!}  2^{2i}
\en
for the number
of walks that make the transition from 0 to 0 in $2i$ steps.
Multiplying by $\cn^{-2i}$ to convert the number of walks into a
probability leads to \refeq{retpr}.  This proves \refeq{Tpij} for
$j=0$, so we take $j \geq 1$.
%
%To prove \refeq{retpr},
%we observe that the left side is equal to the probability that a
%random walk on $\qn$ that starts at the origin returns to the origin
%after $2i$ steps.  This probability is equal to $n^{-2i}$ times the
%number of walks that make the transition from 0 to 0 in $2i$ steps.
%Let $j$ denote the number of dimensions traversed by the random walker,
%so that $j=1, \ldots, i$. Then we can write
%    \eq
%    \lbeq{retpr2}
%    \int |\hat{D}(k)|^{2i} = n^{-2i} \sum_{j=1}^i b_{j,i}(n),
%    \en
%where $b_{j,i}(n)$ denote the number of $2i$-step walks that use
%$j$ dimensions on $\qn$. We claim that
%    \eq
%    \lbeq{retpr3}
%    b_{j,i}(n) \leq n^j c_{j,i}
%    \en
%for some constants $c_{j,i}$ independent of $n$. Indeed, the factor $n^j$
%is an upper bound for the number of possible choices for the $j$ dimensions
%used by the walker. The factors $c_{j,i}$ denote the number of $2i$-step walks in
%$\mathbb{Q}_j$, respectively $\Z^j$, that go from 0 to 0, and this is clearly
%independent of $n$.
%Substituting \refeq{retpr3} into \refeq{retpr2} yields \refeq{retpr}, with
%$a_i=\sum_{j=1}^i c_{j,i}$. Thus, it suffices to show that the second
%integral on the right hand side is bounded uniformly in
%$n\geq 4j+1$. We will do so for $\Z^n$ and $\qn$ separately.

Fix an even integer $s=s(j)$ such that $t=s/(s-1)$ obeys $jt < j + \frac 12$.
By H\"older's inequality,
    \eq
    T_p^\smallsup{i,j}
     \leq
     \left(\int \hat{D}(k)^{is} \right)^{1/s}
     \left( \int \hat{\tau}_p(k)^{jt} \right)^{1/t}.
    \en
By \refeq{retpr}, it suffices to show that $\int \hat{\tau}_p(k)^{jt}$
is bounded by a constant depending on $j$.  We give separate arguments
for this, for $\Z^n$ and $\qn$.

For $\Z^n$, the infrared bound \cite[(4.7)]{HS90a} implies
that $\hat\tau_p(k) \leq 2[1-\hat D(k)]^{-1}$
for sufficiently large $n$, uniformly in $p\leq p_c(\Z^n)$.
Thus,
    \eq
    \lbeq{infj}
    \int \hat{\tau}_p(k)^{jt} \leq 2^{jt}
    \int \frac{1}{[1-\hat{D}(k)]^{jt}} .
    \en
For $A >0$ and $m >0$,
\eq
    \frac{1}{A^m} = \frac{1}{\Gamma(m)}\int_0^\infty u^{m-1} e^{-uA} du,
\en
so that
    \eq
    \lbeq{1Atrick}
    \int  \frac{1}{[1-\hat{D}(k)]^{jt}}
    =
    \frac{1}{\Gamma(jt)}\int_0^\infty du \, u^{jt-1}
    \Big( \int_{-\pi}^\pi e^{-un^{-1}(1-\cos \theta) }
    \frac{d\theta}{2\pi} \Big)^n.
    \en
The right side is non-increasing in $n$, since $\|f\|_p \leq \|f\|_q$
for $0 < p \leq q \leq \infty$ on a probability space.
Since
\eq
    1-\hat{D}(k) = \sum_{j=1}^n (1-\cos k_j) \geq \frac{2}{\pi^2} \frac{|k|^2}{n},
\en
and since $2jt < 2j+1$,
the integral on the left hand side of \refeq{1Atrick} is finite when
$n=2j+1$.
This completes the proof for
$\Z^n$.

For $\qn$, we use the fact that $\hat\tau_p(0) = \chi(p)$ to see that
    \eq
    \lbeq{rwCS}
    \int \hat{\tau}_p(k)^{jt}
    =
    2^{-n}\chi(p)^{jt} + 2^{-n} \sum_{k \in \{0,\pi\}^n : k \neq 0}
    \hat{\tau}_p(k)^{jt}.
    \en
The first term on the right hand side is at most
$2^{-n}\chi(\bar p_c(\qn))^{jt}= 2^{-n}f_n^{jt}$,
which is exponentially small.  For the second term, we recall
from  \cite[Theorem~6.1]{BCHSS04b}
that $\hat\tau_p(k) \leq [1+O(n^{-1})][1-\hat D(k)]^{-1}$, so it
suffices to prove that
    \eq
    \lbeq{sumTpij}
    2^{-n} \sum_{k \in \{0,\pi\}^n : k \neq 0}
    \frac{1}{[1-\hat D(k)]^{jt}}
    \en
is bounded uniformly in $n \geq  1$.

For this, we let $m(k)$ denote the number of nonzero components of $k$.  We
fix an $\varepsilon >0$
and divide the sum according to whether $m(k) \leq \varepsilon n$
or $m(k) >\varepsilon n$.   An elementary computation
(see \cite[Section~2.2.1]{BCHSS04b}) gives
$1-\hat{D}(k) = 2m(k)/n$.
Therefore, the contribution to \refeq{sumTpij} due to $m(k) > \varepsilon n$
is bounded by a constant depending only on $\varepsilon$ and $j$.
On the other hand, for $k\neq 0$, we use
$1-\hat{D}(k) =2m(k)/n \geq 2/n$ to see that
    \eqalign
    2^{-n}\shift\shift
        \sum_{k \in \{0,\pi\}^n:  0<m(k) \leq \varepsilon n}
        \frac{1}{[1-\hat{D}(k)]^{jt}}
    &\leq
    2^{-jt} n^{jt}
    2^{-n} \shift\shift
        \sum_{k \in \{0,\pi\}^n:  0<m(k) \leq \varepsilon n} 1
        \nonumber\\
    &= 2^{-jt} n^{jt} 2^{-n}\sum_{m=1}^{\varepsilon  n} {n \choose m}
    \nnb
\lbeq{nH2}
    & \leq 2^{-jt} n^{jt}  \Pbold (X \leq \varepsilon   n),
    \enalign
where $X$ is a binomial random variable with parameters $(n, 1/2)$.
Since $\Ebold [X] = n/2$, the right side of \refeq{nH2}
is exponentially small in $n$ as $n\rightarrow \infty$
if we choose $\varepsilon < \frac 12$, by standard large
deviation bounds for the binomial distribution
(see, e.g., \cite[Theorem~A.1.1]{AS00}).
This completes the proof for
$\qn$.

\smallskip \noindent {\em Proof of \refeq{Tp}.}
We repeat the argument of \cite[Lemma~5.5]{BCHSS04b} for $\qn$,
which applies verbatim for $\Z^n$.
It follows from the BK inequality that
if $x \neq 0$ then
\eq
\lbeq{taupDtau}
    \tau_p(x) \leq p\cn (D*\tau_p)(x).
\en
Using this, we conclude that
    \eq
    \lbeq{Tppf}
    p\cn(D*\tau_p^{*3})(x)\leq p\cn D(x) + 3(p\cn)^2(D^{*2}*\tau_p^{*3})(x),
    \en
where the first term is the contribution where each of the three two-point
functions $\tau_p(u)$ in $\tau_p^{*3}$ is evaluated at $u=0$, and the second
term takes into account the case
where at least one of the three displacements is nonzero.
Since $p \leq \bar p_c  = \cn^{-1} + O(\cn^{-2}) \leq 2\cn^{-1}$
for large $\cn$, this gives
    \eq
    T_p \leq 2\Omega^{-1} + 12T_p^\smallsup{2,3}
    \leq
    (2+12K_{2,3})\cn^{-1} = K\cn^{-1},
    \en
where in the first inequality we used \refeq{FTinv} to rewrite the second
term of \refeq{Tppf}.

\smallskip \noindent {\em Proof of \refeq{tauiT}.}
For $i \geq 1$, the BK inequality can be applied as in the proof of
\refeq{taupDtau} to obtain
\eq
\lbeq{tauibd}
    \tau_p^\smallsup{i}(x) \leq (p\cn)^i (D^{*i}*\tau_p)(x).
\en
It follows from \refeq{FTinv} and \refeq{tauibd} that
\eq
\lbeq{tauiTpf}
    \sup_x \tau_p^\smallsup{i}(x)
    \leq  \sup_x (p\cn)^i \int \hat{D}(k)^i \hat{\tau}_p(k) e^{-ik\cdot x}
    \leq (p\cn)^i T_p^\smallsup{i,1}
    \leq 2^i K_{i,1}\cn^{-i/2},
\en
where we have used the fact that $p\Omega\leq 2$
for $\cn$ sufficiently large.  For $i=1$, this can be improved by
observing that, for $\cn$ sufficiently large,
\eq
    \tau_p^\smallsup{1}(x) \leq p\cn D(x) + \tau_p^\smallsup{2}(x) \leq
    2\cn^{-1} + 2K_{2,1} \cn^{-1}.
\en
\qed

%%%%%%%%%%%%%%%%%%%%%%%%%%%%%%%%%%%%%%%%%%%%%%%%%%%%%%%%%%%%%
%%%%%%%%%%%%%%%%%%%%%%%%%%%%%%%%%%%%%%%%%%%%%%%%%%%%%%%%%%%%%
\section{Proof of Proposition~\ref{prop-Piasya}}
\label{sec-Pibds}

We now complete the proof of Proposition~\ref{prop-Piasya},
by proving \refeq{Pi0asy}, \refeq{Pi1asy}, \refeq{Pi2asy}
in Sections~\ref{sec-pi0a}, \ref{sec-pi1a}, \ref{sec-pi2a},
respectively.
Throughout this section we fix $p \leq \bar p_c(\gr)$.

%%%%%%%%%%%%%%%%%%%%%%%%%%%%%%%%%%%%%%%%%%%%%%%%%%%%%%%%%%%%%%
%%%%%%%%%%%%%%%%%%%%%%%%%%%%%%%%%%%%%%%%%%%%%%%%%%%%%%%%%%%%%%
\subsection{Expansion for $\hat \Pi_p^{\smallsup{0}}$}
\label{sec-pi0a}

Given a configuration, we say that $x$ is {\em doubly connected
    to}\/ $y$, and we write $x \dbc y$, if $x=y$ or
    if there are at least two bond-disjoint
    paths from $x$ to $y$
    consisting of occupied bonds.
For $\ell \geq 4$,
an $\ell$-\emph{cycle} is a set of bonds that can be written
as $\{\{v_{i-1},v_i\}\}_{1 \leq i \leq \ell}$ with $v_\ell = v_0$ and
otherwise $v_i \neq v_j$ for $i \neq j$, and a \emph{cycle} is
an $\ell$-cycle for some $\ell \geq 4$.
By definition,
%        \eq\lbeq{pi0def}
%        \Pi_p^{\smallsup{0}}(x)=\Pbold_p(0 \dbc  x)
%        - \delta_{0, x},
%        \en
%and
\eq
\lbeq{pi0defz}
    \hat\Pi^\smallsup{0}_p = \sum_{x \neq 0} \Pbold_p (0 \dbc x)
    =
    \sum_{x \neq 0}\Pbold_p ( \exists \; \text{occupied cycle containing $0,x$}).
\en
We decompose the summand into (a) the probability that there exists
an occupied 4-cycle containing $0,x$, plus (b) the probability that there exists
an occupied cycle of length at least $6$ containing $0,x$ and no
occupied 4-cycle containing $0,x$.

The contribution to $\hat\Pi^\smallsup{0}_p$
due to (a) is bounded above by summing $p^4$ over $x\neq 0$ and
over 4-cycles containing $0,x$.  The number
of 4-cycles containing $0$ is $\frac 12 \cn\cn'$, and each
such cycle has three possibilities for $x$.  Therefore
\eq
\lbeq{square1}
    \text{contribution due to (a)} \leq \frac{3}{2} \cn\cn' p^4.
\en
For a lower bound, we apply inclusion-exclusion
and subtract from this upper bound the
sum of $p^7$ over $x\neq 0$ and
over pairs of 4-cycles, each containing $0,x$.  In this case, $x$ must
be a neighbour of $0$, and $p^7$ is the probability of simultaneous
occupation of the two 4-cycles.
There are order $\cn^3$ such pairs of 4-cycles.
Since we already know that $\bar p_c(\gr) \leq O(\cn^{-1})$, this gives
\eq
\lbeq{square2}
    \text{contribution due to (a)} =
    \frac{3}{2} \cn\cn' p^4 + O(\cn^3 p^7)
    = \frac{3}{2} \cn\cn' p^4 + O(\cn^{-4}).
\en

For the contribution due to (b), we use Lemma~\ref{lem-cycle} below.
Given increasing
events $E,F$, we use the standard
notation $E\circ F$ to denote the event
that $E$ and $F$ occur disjointly.  Roughly speaking, $E \circ F$ is
the set of bond configurations for which there exist two disjoint sets
of occupied bonds such that the first set guarantees the occurrence of $E$
and the second guarantees the occurrence of $F$.
The BK inequality
asserts that $\Pbold (E \circ F) \leq \Pbold (E)\Pbold(F)$,
for increasing events $E$ and $F$.
(See \cite[Section~2.3]{Grim99} for a proof, and for
a precise definition of $E \circ F$.)

\begin{lemma}
\label{lem-cycle}
Let $p \leq \bar p_c(\gr)$.
Let $\Pi^\smallsup{0,\ell}_p(x)$ denote the probability that
there is an occupied cycle containing $0,x$, of
length $\ell$ or longer.
Then for $\ell \geq 4$ and for $\cn$ sufficiently large (not depending on $\ell$),
\eq
    \sum_{x \neq 0}\Pi^\smallsup{0,\ell}_p(x)
    \leq (\ell-1) 2^\ell K_{\ell,2} \cn^{-\ell/2}.
\en
\end{lemma}

\proof
Let $\ell \geq 4$, and suppose there exists an occupied cycle containing
$0,x$, of length $\ell$ or longer.
Then there is a $j \in \{1,\ldots, \ell-1\}$ such that
$\{0 \xconn{j} x\} \circ \{ 0 \xconn{\ell-j} x\}$ occurs.
Therefore, by the BK inequality,
\eq
    \Pi^\smallsup{0,\ell}_p(x)
    \leq
    \sum_{j=1}^{\ell-1} \tau^\smallsup{j}_p(x)
    \tau^\smallsup{\ell-j}_p(x).
\en
By \refeq{tauibd}, by the fact that
$p\cn \leq 2$ for $\cn$ sufficiently large, and  by \refeq{Tpij}, it follows that
\eq
    \sum_{x \neq 0}\Pi^\smallsup{0,\ell}_p(x)
    \leq
    (\ell-1) 2^\ell (D^{*\ell} * \tau_p^{*2})(0)
    \leq
    (\ell-1) 2^\ell T_p^\smallsup{\ell,2}
    \leq
    (\ell-1) 2^\ell K_{\ell,2}\cn^{-\ell/2},
\en
as required.
\qed

The contribution due to case (b) is therefore at most
$\sum_{x\neq 0} \Pi_p^\smallsup{0,6}(x)
\leq \bigo(\cn^{-3})$, and hence
\eq
    \hat\Pi^\smallsup{0}_p  =
    \frac{3}{2} \cn\cn' p^4 + O(\cn^{-3}),
\en
which proves \refeq{Pi0asy}.

%%%%%%%%%%%%%%%%%%%%%%%%%%%%%%%%%%%%%%%%%%%%%%%%%%%%%%%%%%%%%%
%%%%%%%%%%%%%%%%%%%%%%%%%%%%%%%%%%%%%%%%%%%%%%%%%%%%%%%%%%%%%%
\subsection{Expansion for $\hat \Pi_p^{\smallsup{1}}$}
\label{sec-pi1a}

To define $\hat \Pi_p^{\smallsup{1}}$,
we need the following definitions.

    \begin{defn}
    \label{def-inon}
    (i)
        Given a bond configuration, vertices $x,y$,
        and a set $A$ of vertices of $\gr$, we
        say $x$ and $y$ are \emph{connected through $A$}, and write
        $x \ct{A} y$, if every
        occupied path connecting $x$ to $y$ has at least one bond
        with an endpoint in $A$.
    \newline
        (ii)
        Given a bond configuration, and a bond $b$, we define
        $\tilde{C}^{b}(x)$ to be the set of vertices connected to $x$
        in the new configuration obtained by setting $b$ to be vacant.
        \newline
    (iii) Given a bond configuration and vertices $x,y$,
    we say that the directed bond $(u,v)$ is {\em pivotal} for $x \conn y$
    if (a) $x \conn y$ occurs when the bond $\{u,v\}$ is set occupied, and
    (b) when $\{u,v\}$ is set vacant $x \conn y$ does not occur,
    but $x\conn u$ and $v \conn y$
    do occur. (Note that there is a distinction
    between the events $\{(u,v)$ is pivotal for $x \conn y\}$
    and $\{(v,u)$ is pivotal for $x \conn y\}
    = \{(u,v)$ is pivotal for $y \conn x\}$.)
    \end{defn}

Let
    \eqalign
    \lbeq{317}
    E'(v, x; A) & = \{ v \ct{A} x\} \cap
    \{\not\exists \; \text{pivotal $(u', v')$ for $v \conn x$ s.t.
    $v \ct{A} u'$} \}.
    \enalign
We will refer to the ``no pivotal'' condition of the
second event on the right hand
side of \refeq{317} as the ``NP'' condition.

By definition,
    \eqalign
    \lbeq{Pi1defa}
    \hat\Pi_p^\smallsup{1}
    & = \sum_{x}p\sum_{(u, v)}
    \Ebold_{\sss 0} \left[ I[0 \dbc u]
    {\mathbb P}_{\sss 1} (E'(v, x;\tilde{C}^{(u,v)}_{\SSS 0}(0))) \right]
    ,
    \enalign
where the sum over $(u,v)$ is a sum over directed bonds.
On the right hand side, the cluster $\tilde{C}^{(u,
v)}_{\SSS 0}(0)$ is random with respect to the expectation
$\Ebold_{\SSS 0}$, so that
$\tilde{C}^{(u,v)}_{\SSS 0}(0)$ should be regarded as a \emph{fixed} set
inside the probability $\Pbold_{\SSS 1}$.
The latter introduces a second percolation model which
depends on the original percolation model via the set
$\tilde{C}^{(u, v)}_{\SSS 0}(0)$.
We use subscripts for
$\tilde{C}$ and the expectations,
to indicate to which expectation $\tilde{C}$ belongs,
and refer to the bond configuration corresponding to expectation~$j$
as the ``level-$j$'' configuration.
We also write $F_j$ to indicate an event $F$ at level-$j$.
Then \refeq{Pi1defa} can be written as
\eqalign
\lbeq{Pi1b}
    \hat\Pi_p^\smallsup{1}
    & = \sum_{x}p\sum_{(u, v)}
    \Pbold^\smallsup{1} \left[ \{0 \dbc u\}_{\sss 0}
    \cap E'(v, x;\tilde{C}^{(u,v)}_{\SSS 0}(0))_{\sss 1} \right],
\enalign
where $\Pbold^\smallsup{1}$ represents the joint expectation of the
percolation models at levels-0 and 1.

We begin with a minor extension
of a standard estimate for $\hat\Pi_p^\smallsup{1}$
(see \cite[Section~4]{BCHSS04b} for related discussion with our present
notation).  Making the abbreviation $\tilde{C}_{\sss 0} =
\tilde{C}^{(u,v)}_{\SSS 0}(0)$,
we may insert within the square brackets on the right hand side of \refeq{Pi1b}
the disjoint union
\eq
    \big(\{u=0\} \cap \{x\in \tilde{C}_{\SSS 0}\}\big)
    \stackrel{\bullet}\bigcup
     \big(\{u=0\} \cap \{x\nin \tilde{C}_{\SSS 0}\}\big)
    \stackrel{\bullet}\bigcup
    \{u \neq 0\}  .
\en
The first term is the leading term and the other two produce error terms.

We first show that the term $\{ u \neq 0\}$ produces an error term.
We define the events
    \eqalign
    F_0(0,u,w,z) &= \{0\conn u\}\circ \{0\conn w\}\circ \{w\conn u\}
    \circ \{w\conn z\},
    \lbeq{Fdefa}\\
    F_1(v,t,z,x) &= \{v\conn t\}\circ \{t\conn z\}
    \circ \{t\conn x\}
    \circ \{z\conn x\}.
    \lbeq{FNdefa}
    \enalign
Note that $F_1(v,t,z,x) = F_0(x,z,t,v)$. Recalling the definition
of $\{x\xconn{j} y\}$ from \refeq{xconnidef}, we also define
   \eqalign
    F_0^\smallsup{j}(0,u,w,z)
    &=
    \bigcup_{j_1+j_2+j_3= j}
    \{0\xconn{j_1} u\}\circ \{0\xconn{j_2} w\}\circ \{w\xconn{j_3} u\}
    \circ \{w\conn z\},
    \lbeq{Fdefaell}\\
    F_1^\smallsup{j}(v,t,z,x) &=
    \bigcup_{j_1+j_2+j_3= j}
    \{v\conn t\}\circ \{t\xconn{j_1} z\}
    \circ \{t\xconn{j_2} x\}
    \circ \{z\xconn{j_3} x\}.
    \lbeq{FNdefaell}
    \enalign

For $u \neq 0$,
it can be seen from the fact that $u$ and $0$ are in a
level-0 cycle of length
at least $4$ that
    \eqalign
    \lbeq{Fbda}
    &\{0 \dbc u \neq 0\}_0 \cap
        E'(v, x; \tilde{C}_{0})_1
    %\\ \nonumber &\quad \quad
    \subset  \bigcup_{t,w,z}
    \Big(
    F_0^\smallsup{4}(0, u, w, z)_0
     \cap
    F_1(v,t,z,x)_1\Big),
    \enalign
and hence this contribution to $\hat\Pi_p^\smallsup{1}$ is at most
    \eqalign
    \lbeq{PiFs}
    p\sum_{x,(u,v),t,w,z}
    \Pro (F_0^\smallsup{4}(0,u,w,z))
    \Pro(F_1(v,t,z,x)).
    \enalign
Let
    \eqalign
        A_3(t,z,x)
        & =
        \tau_{p}(x-t)\tau_{p}(z-t)\tau_{p}(z-x),\lbeq{A3def}\\
        A_3^\smallsup{j}(t,z,x)
        & =
        \sum_{j_1+j_2+j_3 = j}
        \tau_{p}^\smallsup{j_1}(x-t)\tau^\smallsup{j_2}_{p}(z-t)
        \tau^\smallsup{j_3}_{p}(z-x),\lbeq{A3defell}\\
        B_1(w,u,z,t) & =  (p\cn D*\tau_{p})(t-u) \tau_p(z-w).
        \lbeq{B1def}
%       \\
%        B_2(u,v,s,t) & = \tau_{p}(u,v)\tau_{p}(u,t)
%        \tau_{p}(v,s)\tau_{p}(s,t)
%    \nonumber \\ & \quad
%    + \sum_{a\in \Vbold}
%        \tau_{p}(s,a)\tau_{p}(a,u)\tau_{p}(a,t) \delta_{v,s}\tau_{p}(u,t).
%        \lbeq{B2def}
    \enalign
By the BK
inequality, \refeq{PiFs} is at most
    \eqalign
    \lbeq{PibdAB}
    \sum_{u,w}
    A_3^\smallsup{4}(0,u, w)
    \sum_{t,z}
    B_1(w,u,z,t)
    \sum_{x}
    A_3(t,z, x).
    \enalign
Replacing $w,z,t,x$ by $w=w'+u$, $z=z'+u$, $t=t'+u$, $x=x'+u$,
and using symmetry, this
is equal to
   \eqalign
    \lbeq{PibdABa}
    \sum_{u,w'}
    A_3^\smallsup{4}(0,u, w')
    \sum_{t',z'}
    B_1(w',0,z',t')
    \sum_{x'}
    A_3(t',z', x').
    \enalign
We note that $B_1(w',0,z',t')=B_1(-t',-t',z'-w'-t',0)$, and set
$z''=z'-t'$, $x''=x'-t'$ and then $t''=-t'$ to rewrite \refeq{PibdABa} as
   \eqalign
    \lbeq{PibdABb}
    &\sum_{u,w'}
    A_3^\smallsup{4}(0,u, w')
    \sum_{t'',z''}
    B_1(t'',t'',z''-w',0)
    \sum_{x''}
    A_3(0,z'', x'')
    \nnb
    & \quad \leq
    \left( \sum_{u,w'}
    A_3^\smallsup{4}(0,u, w')\right)
    \left(
    \sup_a \sum_{t''}
    B_1(t'',t'',a,0)
    \right)
    \left(
    \sum_{z'',x''}
    A_3(0,z'', x'')
    \right).
    \enalign
By \refeq{Tp-def} and the fact that $\tau_p(u) \leq (\tau_p*\tau_p)(u)$,
\eqalign
    \sup_a \sum_{t''}
    B_1(t'',t'',a,0)
    &= \sup_a ( p\cn D*\tau_p * \tau_p)(a)
    \leq
    T_p .
\enalign
Also,
\eq
    \sum_{z'',x''}
    A_3(0,z'',x'')
    =
    (\tau_p*\tau_p*\tau_p)(0) = T_p^\smallsup{0,3},
\en
and, using \refeq{tauibd} and $p\cn \leq 2$,
\eq
    \sum_{u,w'}
    A_3^\smallsup{4}(0,u,w')
    =
    \sum_{j_1+j_2+j_3 = 4}
        (\tau_{p}^\smallsup{j_1}*\tau^\smallsup{j_2}_{p}*
        \tau^\smallsup{j_3}_{p})(0)
        \leq O(T^\smallsup{4,3}).
\en
Therefore, \refeq{PibdABb} is bounded above by
$O(T^\smallsup{4,3}_p T_p T_p^\smallsup{0,3})$, which is
$O(\cn^{-3})$ by \refeq{Tpij} and \refeq{Tp}.

Similarly, an upper bound $O(\cn^{-3})$ can be obtained for the contribution
due to $\{u=0\} \cap \{x \nin \tilde{C}_{\sss 0}\}$, starting from the
observation that
    \eqalign
    \lbeq{Fbda2}
    &\{u=0\} \cap \{x \nin \tilde{C}_{\sss 0} \} \cap
        E'(v, x; \tilde{C}_{\sss 0})_1
    %\\ \nonumber &\quad \quad
    \subset  \bigcup_{t,w,z: z\neq x}
    \Big(
    F_0 (0, 0, 0, z)_0
     \cap
    F_1^\smallsup{4}(v,t,z,x)_1\Big).
    \enalign
The inclusion \refeq{Fbda2} follows from the fact that
if $ x \nin \tilde{C}_{\SSS 0}$, then
to obtain a non-zero contribution to
${\mathbb P}_{\sss 1} (E'(v, x;\tilde{C}_{\SSS 0}))$,
$x$ must be in a level-1 occupied cycle
of length at least 4 which contains a vertex $z \in \tilde{C}_{\sss 0}$.

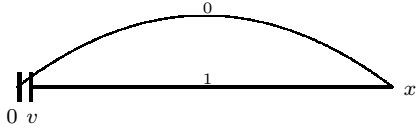
\begin{figure}[t]
\vskip2cm
\begin{center}
\setlength{\unitlength}{0.0075in}%
\begin{picture}(350,0)
\put(43,-5){${\scriptstyle 0}$}
\put(57,-5){${\scriptstyle v}$}
\put(320,15){${\scriptstyle x}$}
\thinlines
%\qbezier(10,20)(40,-10)(70,20)
%\qbezier(10,20)(40,50)(70,20)
\qbezier(50,20)(180,120)(312,20)
%\qbezier(240,5)(360,-50)(480,5)
\Thicklines
\put(52,10){\line(0,1){20}}
\put(60,10){\line(0,1){20}}
\put(60,20){\line(1,0){252}}
%\qbezier(290,20)(320,-10)(350,20)
%\qbezier(290,20)(320,50)(350,20)
%\qbezier(180,20)(375,110)(570,20)
\put(180,22){${\sss 1}$}
%\put(260,22){${\sss 2}$}
%\put(180,5){${\scriptstyle z}$}
\put(180,72){${\sss 0}$}
%\put(375,53){${\sss 3}$}
%\put(445,22){${\sss 4}$}
\thicklines
%\put(325,-5){${\scriptstyle v_{\sss 1}}$}
%\put(305,-5){${\scriptstyle u_{\sss 1}}$}
%\thinlines
%\put(312,10){\line(0,1){20}}
%\put(320,10){\line(0,1){20}}
%\put(320,20){\line(1,0){250}}
%\qbezier(540,20)(570,-10)(600,20)
%\qbezier(540,20)(570,50)(600,20)
\end{picture}
\end{center}
%\vskip1cm

\caption{Depiction of the event appearing in
\refeq{Pi1defux}.  Line~0 corresponds to a connection
in level-0 and line~1 to a connection in level-1.}
\label{fig-1}
\end{figure}

We are left to consider the leading term
    \eqalign
    \lbeq{Pi1defux}
    & \sum_{x}p\sum_{(0, v)}
    \Pbold^{\smallsup{1}}\big(\{x \in\tilde{C}^{(0,v)}_{\SSS 0}(0)\}
    \cap
    E'(v, x;\tilde{C}^{(0,v)}_{\SSS 0}(0))_{\sss 1}\big).
    \enalign
See Figure~\ref{fig-1} for a depiction of the event
appearing in \refeq{Pi1defux}.

%%%%%FIGFIGFIGFIGFIGFIGFIGFIGFIGFIGFIGFIGFIGFIGFIGFIGFIGFIGFIGFIG
%%%%%FIGFIGFIGFIGFIGFIGFIGFIGFIGFIGFIGFIGFIGFIGFIGFIGFIGFIGFIG

The event in \refeq{Pi1defux} is a subset of the event
$\{x\in \tilde{C}_{\sss 0}\} \cap \{v\conn x\}_{\sss 1}$.
Thus, either there is a level-0 connection from $0$ to $x$
(not using the bond $\{0,v\}_{\sss 0}$) of
length $\ell_0$ and a level-1 connection from $v$ to $x$ of length
$\ell_1$, with $\ell_0+\ell_1 \leq 4$, or
$\{0\xconn{i_0} x\}_{\sss 0}\cap \{v\xconn{i_1} x\}_{\sss 1}$
occurs with $i_0+i_1 = 5$.  This decomposition is not disjoint,
as the latter possibility does not
imply that the former does not occur, but this is fine for an upper bound.
By \refeq{tauibd} and \refeq{Tpij}, the contribution due
to the latter case is bounded above by
    \eqalign
    \sum_{i_0=0}^5 \sum_x p\sum_{(0,v)}
    \tau_p^{\smallsup{i_0}}(x)\tau_p^{\smallsup{5-i_0}}(x-v)
    &=
    \sum_{i_0=0}^5
    (p\cn D*\tau_p^{\smallsup{i_0}}*\tau_p^{\smallsup{5-i_0}})(0)
    \nnb &
    \leq
    6(p\cn)^6 T_p^{\smallsup{6,2}}=O(\Omega^{-3}),
    \lbeq{Pi1lellbd}
    \enalign
so this is an error term.

Since $v$ and $0$ have opposite parity,
if there is a level-0 connection from $0$ to $x$ of
length $\ell_0$ and a level-1 connection from $v$ to $x$ of length
$\ell_1$, then $\ell_0+\ell_1$ must be odd.
Thus, we are left to deal with the cases $\ell_0+\ell_1=1$
and $\ell_0+\ell_1=3$, and
we consider these separately.

\smallskip \noindent
{\em Case that $\ell_0+\ell_1=1$.}
If $\ell_1=0$, then $x=v \in \tilde{C}^{(0,v)}_{\sss 0}(0)$, which
forces $\ell_0 \geq 3$.  This is inconsistent with $\ell_0+\ell_1=1$
and therefore need not be considered here.
We may therefore assume that $\ell_0=0$ and $\ell_1=1$,
so that $x=0$,
$\{0,v\}_{\sss 1}$ is occupied, and, to satisfy the NP
condition of \refeq{317},
$v \nin \tilde{C}^{(0,v)}_{\sss 0}(0)$. We use inclusion-exclusion
on the latter, writing
    \eq
    \lbeq{NP1}
    I[v \nin \tilde{C}^{(0,v)}_{\sss 0}(0)]
    =1-I[v \in \tilde{C}^{(0,v)}_{\sss 0}(0)].
    \en
The first term contributes
    \eq
    p\sum_{(0,v)} p = \Omega p^2.
    \en
The second term requires a level-0
connection from $0$ to $v$ of length 3 or more,
which has probability
$\tau_p^\smallsup{3}(v)$, so that by \refeq{tauibd} and \refeq{Tpij},
the second term contributes
    \eq
    p\sum_{(0,v)} p \tau_p^\smallsup{3}(v)
    = p(p\cn)(D*\tau_p^\smallsup{3})(0)
    \leq p(p\cn)^4 T^{\smallsup{4,1}} =O(\Omega^{-3}),
    \en
and hence is an error term. Thus, the case $\ell_0+\ell_1=1$ contributes
    \eq
    \lbeq{contr1}
    \Omega p^2 +O(\Omega^{-3}).
    \en

\smallskip \noindent
{\em Case that $\ell_0+\ell_1=3$.}
There are four possibilities: $\ell_1= 0,1,2 , 3$.
If $\ell_1=
0$ then $x=v$, the NP condition
is trivially satisfied, and there is an occupied level-0 path from 0 to $v$ of
length 3.  This contribution is
    \eq
    \lbeq{contr2}
    \Omega \Omega' p^4+O(\Omega^{3}p^7)
    =\Omega \Omega' p^4+O(\Omega^{-4}),
    \en
where we have used inclusion-exclusion in a manner similar to that
of the argument around
\refeq{square1}--\refeq{square2}.
In more detail, the first term in
\refeq{contr2} accounts for the sum of the probability of an occupied
level-0 path of length 3 from $0$ to $v$, while the second term accounts
for overcounting due to simultaneous occupation of more than one such path.

For $\ell_1= 1,2 , 3$, we note that
\eq
    \{x \in \tilde{C}_{\sss 0}\} \cap E'(v, x;\tilde{C}_{\SSS 0})_{\sss 1}
    = \{x \in \tilde{C}_{\sss 0}\}\cap  \{v \conn x\}_{\sss 1}
    \cap \,{\rm NP} ,
\en
and use
$I[{\rm NP}]=1-I[{\rm NP}^c]$,
to conclude that
    \eq
    \lbeq{NP2}
    I[x \in \tilde{C}_{\sss 0}]
    I[E'(v, x;\tilde{C}_{\SSS 0})_{\sss 1}]
    =I[x \in \tilde{C}_{\SSS 0}]I[\{v \conn x\}_{\sss 1}]
    -I[x\in \tilde{C}_{\SSS 0}]I[\{v \conn x\}_{\sss 1}]I[{\rm NP}^c].
    \en

We first consider the first term on the right hand side of \refeq{NP2}.
In the following, we write $e$ to denote a neighbour of $0$ that is not $\pm v$,
and which will ultimately be summed over.
We again apply
an inclusion-exclusion argument similar to that used for \refeq{contr2},
but do not discuss its details.

The case $\ell_1=1$ corresponds to $\ell_0=2$,
so that $x=v+e$, with the three bonds
$\{0,e\}_{\sss 0}$, $\{e,x\}_{\sss 0}$, $\{x,v\}_{\sss 1}$ each occupied.
This contributes $\cn\cn'p^4$.  Note that in the related configuration
in which $\{0,v\}_{\sss 0}$, $\{v,x\}_{\sss 0}$, $\{x,v\}_{\sss 1}$
are each occupied, the level-0 path $\{0,v\}_{\sss 0}$, $\{v,x\}_{\sss 0}$
from $0$ to $x$ uses the bond $\{0,v\}_{\sss 0}$, and therefore need not
be considered.  For $\Z^n$, the configuration with $x=2v$
and with $\{0,v\}_{\sss 0}$, $\{v,2v\}_{\sss 0}$, $\{v,2v\}_{\sss 1}$
each occupied
need not be considered for the same reason.  (Also, it contributes
$O(\cn p^4) = O(\cn^{-3})$ which is an error term.)

The case $\ell_1=2$ corresponds to $\ell_0=1$,
so that $x=e$, either with the three bonds
$\{0,x\}_{\sss 0}$, $\{x,x+v\}_{\sss 1}$, $\{x+v,v\}_{\sss 1}$ each occupied,
or with the three bonds
$\{0,x\}_{\sss 0}$, $\{0,x\}_{\sss 1}$, $\{0,v\}_{\sss 1}$ each occupied.
This contributes $2\cn\cn'p^4$.
For $\Z^n$, the configuration with $x=-v$ and
with $\{0,-v\}_{\sss 0}$, $\{0,v\}_{\sss 1}$, $\{0,-v\}_{\sss 1}$
each occupied contributes $O(\cn p^4)=O(\cn^{-3})$ and thus is an error term.

The case $\ell_1=3$ corresponds to $\ell_0=0$,
so that $x=0$, with the three bonds
$\{0,e\}_{\sss 1}$, $\{e,e+v\}_{\sss 1}$, $\{e+v,v\}_{\sss 1}$ each occupied.
This contributes $\cn\cn'p^4$.

In summary, the first term on the right
hand side of \refeq{NP2}, with $\ell_1=1,2,3$, contributes
    \eq
    \lbeq{contr3}
    5\Omega \Omega' p^4+O(\Omega^{-3}).
    \en
Next, we consider the effect of the second term in \refeq{NP2},
for $\ell_1=1,2,3$.

For $\ell_1=1$, we have seen above that,
to leading order,
$\{0,e\}_{\sss 0}$, $\{e,x\}_{\sss 0}$, $\{x,v\}_{\sss 1}$
are each occupied.
The only possible pivotal bond for the level-1 connection from $v$
to $x$ is therefore $(v,x)_{\sss 1}$, and thus the
failure of NP requires $v \in \tilde{C}_{\sss 0}$.
This requires a level-0 connection, disjoint from the bonds
$\{0,e\}_{\sss 0}$ and $\{e,x\}_{\sss 0}$, which joins either
$0$ to $v$, $e$ to $v$, or $x$ to $v$.  This adds an additional
factor $O(\cn^{-1})$ and hence produces an error term.

For $\ell_1=2$, we have seen above that there are two cases to consider.
Suppose first that
$\{0,x=e\}_{\sss 0}$, $\{x,x+v\}_{\sss 1}$, $\{x+v,v\}_{\sss 1}$
are each occupied.  The only possible pivotal bonds for the level-1
connection from $v$ to $x$ are $(v,x+v)_{\sss 1}$
and $(x+v,x)_{\sss 1}$.  Violation of NP therefore requires either
$(v,x+v)_{\sss 1}$ is pivotal and $v \in \tilde{C}_{\sss 0}$,
or $(x+v,v)_{\sss 1}$ is pivotal and $x+v \in \tilde{C}_{\sss 0}$.
In either of these cases, the condition that $\tilde{C}_{\sss 0}$
contain an additional vertex is a higher order effect and leads
to an error term $O(\cn^{-3})$.

The remaining case for $\ell_1=2$ has
$\{0,x=e\}_{\sss 0}$, $\{0,x\}_{\sss 1}$, $\{0,v\}_{\sss 1}$ each occupied.
The only possible pivotal bonds for the level-1
connection from $v$ to $x$ are $(v,0)_{\sss 1}$
and $(0,x)_{\sss 1}$.
Violation of NP therefore requires either
$(v,0)_{\sss 1}$ is pivotal and $v \in \tilde{C}_{\sss 0}$,
or $(0,x)_{\sss 1}$ is pivotal and $0 \in \tilde{C}_{\sss 0}$.
The first of these cases leads to an error term as above.
For the second case, $0 \in \tilde{C}_{\sss 0}$ is automatic,
and inclusion-exclusion applied to the requirement that
$(0,v)_{\sss 1}$ is pivotal leads to a net contribution for
$\ell_1=2$ of  $-\cn\cn'p^4 +O(\cn^{-3})$.

Finally, we consider $\ell_1=3$.
In this case, $x=0$, and
$\{0,e\}_{\sss 1}$, $\{e,e+v\}_{\sss 1}$, $\{e+v,v\}_{\sss 1}$ are each occupied.
The only possible violations of NP are:
$(v,v+e)_{\sss 1}$ is pivotal for the connection from $v$ to $x=0$
and $v \in \tilde{C}_{\sss 0}$, or
$(v+e,e)_{\sss 1}$ is pivotal
and $v +e \in \tilde{C}_{\sss 0}$, or
$(e,0)_{\sss 1}$ is pivotal
and $e \in \tilde{C}_{\sss 0}$.  In any of these three cases,
the condition that $\tilde{C}_{\sss 0}$ must contain the additional
vertex requires extra connections that produce an error term
$O(\cn^{-3})$ overall, using reasoning analogous to that employed
above.

%We first fix the length of line 1 to be $i$, where
%$i\in\{ 1,2, 3\}$. Thus, the length of line 0 is $3-i$. Thus, there must be
%a path of occupied bonds $b_1, b_2, b_3$ such that $b_i=(s,x)$ and $\{b_1, \ldots, b_i
%\text{ occ.}\}_{\sss 1}\cap \{b_{i+1}, \ldots, b_3\text{ occ.}\}_{\sss 0}$ holds.
%When the NP-condition fails, then there must be
%paths $\omega_0$ and $\omega_1$ on levels 0, respectively 1, that are bond disjoint from
%$\{b_1, \ldots, b_i\}$, respectively, $\{b_{i+1}, \ldots, b_3\text{ occ.}\}_{\sss 0}$,
%such that (i) the end points of $\omega_0$ and $\omega_1$ agree; (ii) $\omega_0$ or $\omega_1$
%does not shrink to a point; (iii) the starting point of $\omega_0$ is one of the vertices
%in the bonds $b_{i+1}, \ldots, b_3$ and the starting point of $\omega_1$ is one of the vertices
%in the bonds $b_1, \ldots, b_i$. Thus, by the BK inequality, we can extract a factor
%$\sum_z \tau^{\smallsup{1}}_p(z-s)\tau_p(z-t)$, where $s,t$ are any of the vertices of
%$b_1, b_2$ or $b_3$. By \refeq{tauiT}, the contribution due to $z=t$
%is bounded by $O(\Omega^{-1})$ uniformly in $s,t$. For the contribution due to $z=t$,
%we use that $\tau_p(z-t)=\tau_p^{\smallsup{1}}(z-t),$
%together with \refeq{tauibd} to arrive at
%    \eq
%    \sup_{s,t} \sum_{z}\tau^{\smallsup{1}}_p(z-s)\tau^{\smallsup{1}}_p(z-t)
%    \leq (p\Omega)^2 T^{\smallsup{2,2}} =O(\Omega^{-1}).
%    \en
%Since the number of possible $s,t$ is bounded by $6$, we obtain a bound of the form
%    \eq
%    \lbeq{contr4}
%    p^4 \Omega^2 O(\Omega^{-1})=O(\Omega^{-3}).
%    \en

We have thus shown that the case $\ell_0+\ell_1=3$ yields a net contribution
\eq
\lbeq{sumis3}
    5\cn\cn'p^4 -\cn\cn'p^4 + O(\cn^{-3})
    =
    4\cn\cn'p^4 + O(\cn^{-3}).
\en

\medskip
In summary, combining \refeq{contr1} and \refeq{sumis3},
we have proved \refeq{Pi1asy}, namely
\eq
    \hat\Pi^\smallsup{1}_p = \cn p^2 + 4 \cn \cn' p^4 +O(\cn^{-3}).
\en

\subsection{Expansion for $\hat \Pi_p^{\smallsup{2}}$}
\label{sec-pi2a}

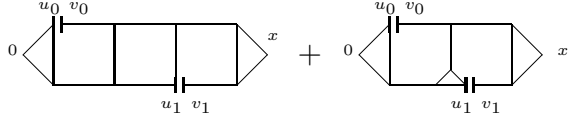
\begin{figure}
\begin{center}
\setlength{\unitlength}{0.008in}%
\begin{picture}(400,90)(20,635)
\thicklines
%\put(105,765){\line( 0,-1){ 10}}
%\put(100,765){\line( 0,-1){ 10}}
%\thinlines
%\put(105,760){\line( 1, 0){ 35}}
%\put(100,720){\line( 1, 0){ 40}}
\thicklines
\put(100,685){\line( 0,-1){ 10}}
\put(105,685){\line( 0,-1){ 10}}
\thinlines
\put(105,680){\line( 1, 0){115}}
\put(140,680){\line( 0,-1){ 40}}
\put(180,680){\line( 0,-1){ 40}}
\thicklines
\put(180,645){\line( 0,-1){ 10}}
\put(185,645){\line( 0,-1){ 10}}
\thinlines
\put(185,640){\line( 1, 0){ 35}}
\put(100,640){\line( 1, 0){ 80}}
\put(320,640){\line( 1, 0){ 30}}
\put(360,680){\line( 0,-1){ 30}}
\thicklines
\put(320,685){\line( 0,-1){ 10}}
\put(325,685){\line( 0,-1){ 10}}
%\put(325,675){\line( 0, 1){  5}}
\thinlines
\put(325,680){\line( 1, 0){ 75}}
\put(370,645){\line( 0, 1){  0}}
\thicklines
\put(370,645){\line( 0,-1){ 10}}
\put(375,645){\line( 0,-1){ 10}}
\thinlines
\put(375,640){\line( 1, 0){ 25}}
%\put(140,720){\line( 1, 1){ 20}}
%\put(160,740){\line( 0, 1){  0}}
%\put(160,740){\line(-1, 1){ 20}}
%\put(140,760){\line( 0,-1){ 40}}
%\put(100,760){\line( 0,-1){ 40}}
%\put(100,720){\line(-1, 1){ 20}}
%\put( 80,740){\line( 1, 1){ 20}}
\put(100,680){\line( 0,-1){ 40}}
\put(100,640){\line(-1, 1){ 20}}
\put( 80,660){\line( 1, 1){ 20}}
\put(220,680){\line( 0,-1){ 40}}
\put(220,640){\line( 1, 1){ 20}}
\put(240,660){\line(-1, 1){ 20}}
%\multiput(360,650)(-0.40000,-0.40000){26}{\makebox(0.4444,0.6667){.}}
\put(360,650){\line(-1,-1){ 10}}
\put(350,640){\line( 1, 0){ 20}}
%\multiput(370,640)(-0.40000,0.40000){26}{\makebox(0.4444,0.6667){\sevrm .}}
\put(360,650){\line( 1,-1){ 10}}
\put(400,680){\line( 0,-1){ 40}}
\put(400,640){\line( 1, 1){ 20}}
\put(420,660){\line(-1, 1){ 20}}
\thinlines
\put(320,680){\line( 0,-1){ 40}}
\put(320,640){\line(-1, 1){ 20}}
\put(300,660){\line( 1, 1){ 20}}
\thinlines
\put( 70,660){\makebox(0,0)[lb]{\raisebox{0pt}[0pt][0pt]{${\sss 0}$}}}
\put(290,660){\makebox(0,0)[lb]{\raisebox{0pt}[0pt][0pt]{${\sss 0}$}}}
\put(260,655){\makebox(0,0)[lb]{\raisebox{0pt}[0pt][0pt]{$+$}}}
\put( 240,670){\makebox(0,0)[lb]{\raisebox{0pt}[0pt][0pt]{${\sss x}$}}}
\put(430,660){\makebox(0,0)[lb]{\raisebox{0pt}[0pt][0pt]{${\sss x}$}}}
\put( 90,690){\makebox(0,0)[lb]{\raisebox{0pt}[0pt][0pt]{${\sss u_0}$}}}
\put(310,690){\makebox(0,0)[lb]{\raisebox{0pt}[0pt][0pt]{${\sss u_0}$}}}
\put( 110,690){\makebox(0,0)[lb]{\raisebox{0pt}[0pt][0pt]{${\sss v_0}$}}}
\put(330,690){\makebox(0,0)[lb]{\raisebox{0pt}[0pt][0pt]{${\sss v_0}$}}}
\put( 170,625){\makebox(0,0)[lb]{\raisebox{0pt}[0pt][0pt]{${\sss u_1}$}}}
\put(360,625){\makebox(0,0)[lb]{\raisebox{0pt}[0pt][0pt]{${\sss u_1}$}}}
\put( 190,625){\makebox(0,0)[lb]{\raisebox{0pt}[0pt][0pt]{${\sss v_1}$}}}
\put(380,625){\makebox(0,0)[lb]{\raisebox{0pt}[0pt][0pt]{${\sss v_1}$}}}
\end{picture}
\end{center}
\caption{\label{fig-pi2diag} The standard
diagrams bounding $\Pi^\smallsup{2}_p$.  All vertices other than $0$ are
summed over the vertex set $\mathbb V$ of $\gr$, lines represent factors
of $\tau_p$, and vertical bars represent factors
$p\cn D$. }
\end{figure}

By definition,
    \eq
    \lbeq{Pi2def}
    \hat\Pi_p^\smallsup{2}=\sum_{x}\sum_{(u_0, v_0)}\sum_{(u_1, v_1)}
    p^2\Ebold_{\sss 0} \Big[I[0 \dbc u_0]
    {\mathbb E}_{\sss 1} \big[I[E'(v_0, u_1;\tilde{C}_{\SSS 0})]
    {\mathbb E}_{\sss 2} I[E'(v_1, x;\tilde{C}_{\SSS 1})]
    \big]\Big],
    \en
where we have made the abbreviations $\tilde{C}_{\SSS 0}
=\tilde{C}^{(u_0,v_0)}_{\SSS 0}(0)$ and
$\tilde{C}_{\SSS 1} =\tilde{C}^{(u_1,v_1)}_{\SSS 1}(v_0)$.
A standard estimate for $\hat\Pi_p^\smallsup{2}$ is
    \eq
    \lbeq{Pi2Tbd}
    0 \leq \hat\Pi_p^\smallsup{2}
    \leq 2T_p^\smallsup{0,3} (T_p T_p^\smallsup{0,3} )^2
    \en
(see, e.g., \cite[Section~4.2]{BCHSS04b}; one factor 2 in
\cite[Proposition~4.1]{BCHSS04b} is easily dropped for $N=2$).
This estimate arises from the upper bound for $\hat\Pi^\smallsup{2}_p$
depicted in Figure~\ref{fig-pi2diag}.  The factor 2 is due to the
fact that there are two terms in the upper bound.  The two factors
$T_p$ in each term
arise from the two diagram loops containing lines with vertical
bars, and the three factors $T_p^\smallsup{0,3}$ arise from the
other three diagram loops.

We claim that contributions to $\hat\Pi^\smallsup{2}_p$ in which
$u_0 \neq 0$, or $u_1 \nin \tilde{C}_{\sss 0}$, or $x \nin \tilde{C}_{\sss 1}$
produce an error term of order $O(\cn^{-4})$.  This follows from
routine estimates, along the lines of those used in Section~\ref{sec-pi1a}
to conclude that we could assume there that $u=0$ and $x \in \tilde{C}_{\sss 0}$.
These estimates, which we do not write down here in detail, show
for example that if
$u_0\neq 0$, then the factor $T_p^\smallsup{0,3}$,
arising from the leftmost diagram loop can be replaced by a constant
multiple of $T_p^\smallsup{4,3}.$  By \refeq{Tpij} and \refeq{Tp},
this leads to a bound $O(\Omega^{-4})$, which is an error term.
Similarly, if $u_1 \nin \tilde{C}^{(0,v_0)}_{\SSS 0}(0)$, then
the event
$E'(v_0, u_1;\tilde{C}^{(0,v_0)}_{\SSS 0}(0))$ requires that
$u_1$ must be in an occupied level-1 cycle of length at least 4.
In this case, we may again use standard
estimates to replace a factor
$T_p^\smallsup{0,3}$ in \refeq{Pi2Tbd}, arising from the diagram loop in
Figure~\ref{fig-pi2diag} containing $u_1$,
by a constant multiple of $T_p^\smallsup{4,3}$,
and again this contribution is $O(\cn^{-4})$. Finally, the same
situation arises when $x\nin \tilde{C}^{(u_1,v_1)}_{\SSS 0}(v_0)$,
in which case we can replace  the factor
$T_p^\smallsup{0,3}$ arising from the rightmost
diagram loop by a constant multiple of $T_p^\smallsup{4,3}$,
and again this contribution is $O(\cn^{-4})$. Thus, we are now
left to analyze
    \eq
    \lbeq{Pi2defa}
    \sum_{x}\sum_{(0, v_0)}\sum_{(u_1, v_1)}
    p^2\Pbold^{\smallsup{2}} \Big(\{u_1\in\tilde{C}_{\SSS
    0}\}\cap \{x\in \tilde{C}_{\SSS 1}\}\cap
    E'(v_0, u_1;\tilde{C}_{\SSS 0})_{\sss 1}\cap
    E'(v_1, x;\tilde{C}_{\SSS 1})_{\sss 2}\Big),
    \en
where we write $\Pbold^{\smallsup{2}}$ for the joint probability
of levels 0, 1 and 2. Let $G$ denote the
intersection of events on the right hand
side of \refeq{Pi2defa}.

%%%%%FIGFIGFIGFIGFIGFIGFIGFIGFIGFIGFIGFIGFIGFIGFIGFIGFIGFIGFIGFIG
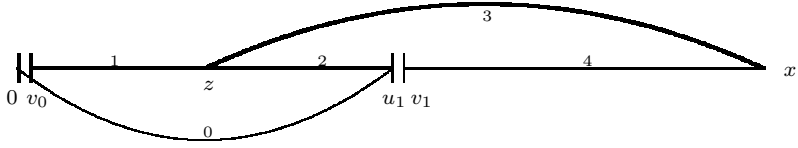
\begin{figure}[t]
\vskip2cm
\begin{center}
\setlength{\unitlength}{0.0075in}%
\begin{picture}(700,-100)
\put(43,-5){${\scriptstyle 0}$}
\put(57,-5){${\scriptstyle v_{\sss 0}}$}
\put(585,15){${\scriptstyle x}$}
\thinlines
%\qbezier(10,20)(40,-10)(70,20)
%\qbezier(10,20)(40,50)(70,20)
\qbezier(50,20)(180,-80)(312,20)
%\qbezier(240,5)(360,-50)(480,5)
\Thicklines
\put(52,10){\line(0,1){20}}
\put(60,10){\line(0,1){20}}
\put(60,20){\line(1,0){252}}
%\qbezier(290,20)(320,-10)(350,20)
%\qbezier(290,20)(320,50)(350,20)
\qbezier(180,20)(375,110)(570,20)
\put(115,22){${\sss 1}$}
\put(260,22){${\sss 2}$}
\put(180,5){${\scriptstyle z}$}
\put(180,-28){${\sss 0}$}
\put(375,53){${\sss 3}$}
\put(445,22){${\sss 4}$}
\thicklines
\put(325,-5){${\scriptstyle v_{\sss 1}}$}
\put(305,-5){${\scriptstyle u_{\sss 1}}$}
\thinlines
\put(312,10){\line(0,1){20}}
\put(320,10){\line(0,1){20}}
\put(320,20){\line(1,0){250}}
%\qbezier(540,20)(570,-10)(600,20)
%\qbezier(540,20)(570,50)(600,20)
\end{picture}
\end{center}
\vskip1cm

\caption{Diagrammatic representation of the event
\refeq{disjconn}. Line~0 corresponds to a connection
in level-0, lines~1, 2, 3 correspond to connections in level-1 and
line~4 to a connection in level-2.}
\label{fig-2}
\end{figure}
%%%%%FIGFIGFIGFIGFIGFIGFIGFIGFIGFIGFIGFIGFIGFIGFIGFIGFIGFIGFIG

The event $G$ on the right hand
side of \refeq{Pi2defa} is contained in the event
    \eq
    \lbeq{disjconn}
    \{0\conn u_1\}_{\sss 0}\cap \bigcup_{z\in \gr}
    \big\{\{v_0\conn z\}\circ \{z\conn
    v_1\}\circ\{z\conn x\}\big\}_{\sss 1} \cap \{v_1\conn
    x\}_{\sss 2},
    \en
which is depicted in Figure~\ref{fig-2}.
For any choice of $j_0,j_1,j_2,j_3,j_4$,
a subset of \refeq{disjconn} is the event
    \eq
    \lbeq{disjconn2}
    \{0\xconn{j_0} u_1\}_{\sss 0}\cap \bigcup_{z\in \gr}
    \big\{\{v_0\xconn{j_1} z\}\circ \{z\xconn{j_2}
    v_1\}\circ\{z\xconn{j_3} x\}\big\}_{\sss 1} \cap \{v_1\xconn{j_4}
    x\}_{\sss 2}.
    \en
Since $v_0$ has odd parity, and since $u_1$ and $v_1$ have
opposite parity, we may assume that $j_0+j_1+j_2$ and $j_2+j_3+j_4$
are both odd.
If $j_0+j_1+j_2 \geq 3$, then a standard diagrammatic
estimate gives $O(T^\smallsup{4,3}_pT_p) = O(\cn^{-3})$ for the
contribution of \refeq{disjconn2}.
Similarly, if
$j_2+j_3+j_4 \geq 3$, then again a standard diagrammatic estimate
gives an upper bound $O(T_pT^\smallsup{4,3}_p)=O(\cn^{-3})$.
Note that if $u_1 \neq 0$, then we may assume that $j_0+j_1+j_2 \geq 3$,
which gives an error term.

%so that  we obtain
%    \eqalign
%    \refeq{Pi2defa} &\leq \sum_{(0,v_0), (u_1,v_1), x, z}
%    p^2
%    \tau_p(u_1)\tau_p(z-v_0)\tau_p(v_1-z)\tau_p(x-z)\tau_p(x-v_1)\nn\\
%    &\leq \sum_{v_0, (u_1,v_1), z}
%    p^2
%    \tau_p(u_1)\tau_p(z-v_0)\tau_p(v_1-z)\sup_{z-v_1}\sum_x
%    \tau_p(x-z)\tau_p(x-v_1)\nn\\
%    &\leq T_p^{\smallsup{1,3}}T_p^{\smallsup{1,2}},
%    \lbeq{Pi2defabd}
%    \enalign
%which is indeed a refinement of \refeq{Pi2Tbd}.
%We will now refine this bound further by
%bounding the contribution of the event
%By the BK inequality and the independence of
%levels 0, 1, 2, the contribution of \refeq{disjconn2} to
%\refeq{Pi2defa} is bounded above by
%    \eqalign
%    &\sum_{x, z} \sum_{(0,v_0)}\sum_{(u_1,v_1)}
%    p^2
%    \tau_p^{\smallsup{j_0}}(u_1)\tau_p^{\smallsup{j_1}}(z-v_0)
%    \tau_p^{\smallsup{j_2}}(v_1-z)\tau_p^{\smallsup{j_3}}(x-z)
%    \tau_p^{\smallsup{j_4}}(x-v_1)
%    \enalign
%Using \refeq{tauibd} and
%Lemma~\ref{lem-Tbd}, this can be bounded above by
%\eqalign
%    & T_p^{\smallsup{1+j_0+j_1+j_2,3}}
%    T_p^{\smallsup{1+j_2+j_3,2}}
%    =O(\cn^{-3}),
%    \lbeq{Pi2defabdell}
%    \enalign
%which is an error term.

Thus, we may assume that $G$ occurs, that $u_1=0$,
and that there is a $z$ such
that the connections of Figure~\ref{fig-2} occur with lines of
length $\ell_0=0,\ell_1,\ell_2,\ell_3,\ell_4$, where
\eqalign
     \ell_1+\ell_2 & = 1,
    \\
     \ell_2+\ell_3+\ell_4 & = 1.
\enalign
This gives three possibilities for $(\ell_0,\ell_1,\ell_2,\ell_3,\ell_4)$,
namely
\eq
    (0,0,1,0,0), \quad (0,1,0,1,0), \quad (0,1,0,0,1),
\en
and it suffices to compute the contribution from each of these cases.

\smallskip \noindent
{\em Case of $(0,0,1,0,0)$.}
In this case, $u_0=u_1=0$, $z=x=v_1=v_0$, and the bond
$\{v_0,u_1\}_{\sss 1}$ is occupied.
We examine the constraints imposed by the event $G$ of \refeq{Pi2defa}.
The events $\{ u_1 \in\tilde{C}_{\sss 0}\}$ and
$\{ x \in\tilde{C}_{\sss 1}\}$ occur trivially.
For the event
$E'(v_0,u_1;\tilde{C}_{\sss 0})_{\sss 1}$, we note that
$\{v_0 \ct{\tilde{C}_{\sss 0}} u_1\}_{\sss 1}
=\{v_0 \ct{\tilde{C}_{\sss 0}} 0\}_{\sss 1}$ occurs.
Violation of the
NP condition requires $\{0 \xconn{3} v_0\}_{\sss 0}$, and this contributes
at most $\sum_{(0,v_0)} p^3 \tau^\smallsup{3}_p(v_0)
\leq p^2 (p\cn)^4 T^\smallsup{4,1}\leq O(\cn^{-4})$.
Thus, up to an error term, we may assume that
$E'(v_0,u_1;\tilde{C}_{\sss 0})_{\sss 1}$ occurs.
Finally, the event
$E'(v_1, x;\tilde{C}_{\SSS 1})_{\sss 2}$ occurs trivially,
since $v_1=x \in \tilde{C}_{\SSS 1}$.
This case contributes
    \eq
    \cn p^3 +\bigo(\cn^{-3}).
    \lbeq{2-a}
    \en

\smallskip \noindent
{\em Case of $(0,1,0,1,0)$.}
In this case, $u_1=0$, $z=u_1$, $v_1=x$.  Also, the fact that
$x \in \tilde{C}_{\sss 1}$ implies that there must be an
occupied level-1 path from $v_0$ to $z=u_1$ to $x=v_1$ that does not
use the bond $(u_1,v_1)_{\sss 1}$.  This implies that
the event $\{ u_1 \xconn{3} v_1\}_{\sss 1}$ occurs, and hence this
case contributes an error term because it corresponds to \refeq{disjconn2}
with $j_3=3$.

\smallskip \noindent
{\em Case of $(0,1,0,0,1)$.}
In this case, $x=z=u_1=u_0=0$, and the bonds $\{0,v_0\}_{\sss 1}$,
$\{0,v_1\}_{\sss 2}$ are occupied.
We denote the neighbours of
$0$ by $e_l$ ($l=1,\ldots, \cn$), so that $v_0=e_i$ and $v_1=e_j$
for some $i,j$.
We examine the constraints imposed by the event $G$ of \refeq{Pi2defa}.
The event $\{u_1 \in \tilde{C}_{\sss 0}\}$
is satisfied trivially, since $u_1=0$.
For the event $\{x \in \tilde{C}_{\sss 1}\}$, we consider separately
the cases $i=j$ (i.e., $v_0=v_1$) and $i \neq j$ (i.e., $v_0 \neq v_1$).
If $i=j$, then $\{x \in \tilde{C}_{\sss 1}\}$ requires
that $\{v_0 \xconn{3} x\}_{\sss 1}$, so this is an error term
in which \refeq{disjconn2} occurs with $j_1+j_2+j_3 \geq 3$ and $j_4=1$
(in more detail, these inequalities imply either that $j_1=2$, in which
case $j_0+j_1+j_2\geq 3$ since the sum must be odd, or that $j_1 \leq 1$,
which implies that $j_2+j_3+j_4 \geq 3$).
If $i\neq j$, then $\{x \in \tilde{C}_{\sss 1}\}$ is achieved by
the bond $\{x,v_0\}_{\sss 1}=\{0,v_0\}_{\sss 1}$.  Thus, we assume
henceforth that $i \neq j$.

For the $E'$ events, we first note that
$\{v_0 \ct{\tilde{C}_{\sss 0}} u_1\}_{\sss 1}$ occurs, since $u_1=0$,
$\{0,v_0\}_{\sss 1}$ is occupied, and $0 \in \tilde{C}_{\sss 0}$.
Also, $\{v_1 \ct{\tilde{C}_{\sss 0}} x\}_{\sss 2}$ occurs, since $x=0$,
$\{0,v_1\}_{\sss 2}$ is occupied, and $0 \in \tilde{C}_{\sss 1}$
(when $i \neq j$).  We will argue below that the NP condition in each
$E'$ event can be neglected, up to an error term.  Assuming this,
this case contributes
\eq
    p\sum_{(0,v_0)} p\sum_{(0,v_1): \, v_1\neq v_0} p^2 +O(\cn^{-3})
    =
    \cn(\cn-1)p^4+\bigo(\cn^{-3}).
    \lbeq{2-b}
\en

If the NP condition is violated for $E'(v_0,u_1;\tilde{C}_{\sss 0})_{\sss 1}
=E'(v_0,0;\tilde{C}_{\sss 0})_{\sss 1}$,
then the bond $(v_0,0)_{\sss 1}$ must be pivotal for the level-1 connection
from $v_0$ to $0$, and moreover $v_0 \in \tilde{C}_{\sss 0}
=\tilde{C}_{\sss 0}^{(0,v_0)}(0)$ must occur.
The latter gives an additional factor $\tau_p^\smallsup{3}(v_0)
\leq O(\cn^{-3/2})$, and hence this contributes to an error term.

If the NP condition is violated for $E'(v_1,x=0;\tilde{C}_{\sss 1})_{\sss 2}$,
then the bond $(v_1,0)_{\sss 2}$ must be pivotal for the level-2 connection
from $v_1$ to $0$, and also $v_1 \in \tilde{C}_{\sss 1}
=\tilde{C}_{\sss 1}^{(0,v_1)}(v_0)$ must occur.  The latter
gives an additional factor $\tau_p^\smallsup{2}(v_1-v_0)
\leq O(\cn^{-1})$, and hence this contributes to an error term.

\medskip
Combining \refeq{2-a}--\refeq{2-b}, we have
\eq
    \hat\Pi_p^\smallsup{2} = \cn p^3 + \cn(\cn-1)p^4+\bigo(\cn^{-3}),
\en
which is \refeq{Pi2asy}.

\section{Conclusions}

We have used the lace expansion to prove that
$p_c(\gr) = \cn^{-1}+\cn^{-2}+\frac 7 2 \cn^{-3} +O(\cn^{-4})$
for $\gr = \Z^n$ and $\gr = \qn$.  This extends by two terms the result
$p_c(\qn)= n^{-1} +O(n^{-2})$ of \cite{BCHSS04b},
and gives a simplified proof of a result of
\cite{HS93up,HS95} for $\Z^n$.

Our proof is essentially mechanical, and with sufficient labour
could be directly extended to compute higher coefficients.  In particular,
it would be
interesting to compute the coefficient of $\cn^{-4}$, which we
expect will be different for $\Z^n$ and $\qn$.

We expect that our method can also be applied to other finite
graphs for which the lace expansion has been proved
to converge in \cite{BCHSS04b}.  A specific example is
the Hamming cube, which has vertex set $\{0,1,\ldots,s\}^n$ with
$s \geq 1$ fixed, and edge set consisting of pairs of vertices which
differ in exactly one component.
For $s = 1$, the Hamming cube is the $n$-cube.  For $s \geq 2$,
the Hamming cube contains cycles of
length 3 (in contrast to $\Z^n$ and $\qn$),
and it would be interesting to study their effect
on the expansion coefficients.

\section*{Acknowledgements}
We thank Christian Borgs, Jennifer Chayes and Joel Spencer for
many stimulating discussions related to this work. The work of
RvdH was supported in part by Netherlands Organisation for
Scientific Research (NWO), and was carried out in part at Delft
University of Technology, at the University of British Columbia,
and at Microsoft Research. The work of GS was supported in part by
NSERC of Canada, by a Senior Visiting Fellowship at the Isaac
Newton Institute funded by EPSRC Grant N09176, by EURANDOM, and by
the Thomas Stieltjes Institute.

%\bibliography{../bibdef/bib}
%\bibliography{bib}

\end{document}

%% file: def99c.tex
\def\UseSection{%%
        \numberwithin{equation}{section}
    \theoremstyle{plain}% default theorem style 
        \newtheorem{theorem}    {Theorem}[section]
        \DefineTheorems % Use this to define other environments to be 
}

\def\DefineTheorems{%%
    
    \newtheorem{lemma}      [theorem] {Lemma}
    
    \newtheorem{prop}       [theorem] {Proposition}
    
    \newtheorem{cor}        [theorem] {Corollary}

    \theoremstyle{definition}% ``defn'' theorem style 
    \newtheorem{defn}       [theorem] {Definition}

    \theoremstyle{definition}% ``remark'' theorem style 

}

\newcommand{\bt}   {\begin{theorem}}
\newcommand{\et}   {\end  {theorem}}
\newcommand{\bl}   {\begin{lemma}}
\newcommand{\el}   {\end  {lemma}}
\newcommand{\bp}   {\begin{prop}}
\newcommand{\ep}   {\end  {prop}}
\newcommand{\bc}   {\begin{cor}}
\newcommand{\ec}   {\end  {cor}}
\newcommand{\bd}   {\begin{defn}}
\newcommand{\ed}   {\end  {defn}}

\newcommand{\ba}   {\begin{array}}
\newcommand{\ea}   {\end  {array}}
\newcommand{\be}   {\begin{enumerate}}
\newcommand{\ee}   {\end  {enumerate}}
\newcommand{\bi}   {\begin{itemize}}
\newcommand{\ei}   {\end  {itemize}}

\def\eq#1\en{\begin{equation}#1\end{equation}}  
\def\eqsplit#1\ensplit{
    \begin{equation}\begin{split}#1\end{split}\end{equation}
    }
\def\eqalign#1\enalign{
    \begin{align}#1\end{align}
    }
\def\eqmul#1\enmul{
    \begin{multline}#1\end{multline}
    }
\newcommand{\eqarrstar} {\begin{eqnarray*}} 
\newcommand{\enarrstar} {\end{eqnarray*}} 
\newcommand{\eqarray}   {\begin{eqnarray}} 
\newcommand{\enarray}   {\end{eqnarray}} 
\newcommand{\nnb}   {\nonumber \\} 

\newcommand{\lbeq}[1]  {\label{e:#1}}
\newcommand{\refeq}[1] {\eqref{e:#1}}    % AMS-LaTeX trick!
\makeatletter
\newcommand{\labelcounter}[2]{{%
    \stepcounter{#1}%   First, increase the ``countC'' by one.
    \protected@write\@auxout{}%
    {\string\newlabel{#2}{{\csname the#1\endcsname}{\thepage}}}%
    {\ref{#2}}% Finally, make sure to refer to this label, 
    }}
\makeatother
\newcommand{\sss}   { \scriptscriptstyle } 

\newcommand{\Ebold} {{\mathbb E}}

\newcommand{\Pbold} {{\mathbb P}}

\newcommand{\Vbold} {{\mathbb V}}

\newcommand{\Zbold} {{\mathbb Z}}

\newcommand{\spose}[1] {{\hbox to 0pt{#1\hss}} }
\newcommand{\ltapprox} {\mathrel{\spose{\lower 3pt\hbox{$\mathchar"218$}}
 \raise 2.0pt\hbox{$\mathchar"13C$}}}
\newcommand{\gtapprox} {\mathrel{\spose{\lower 3pt\hbox{$\mathchar"218$}}
 \raise 2.0pt\hbox{$\mathchar"13E$}}}

\newcommand{\nin}  {{ \not\in }}

%% file: ncube4fin.bbl
\begin{thebibliography}{10}
\bibliographystyle{plain}

\bibitem{AB87}
M.~Aizenman and D.J. Barsky.
\newblock Sharpness of the phase transition in percolation models.
\newblock {\em Commun. Math. Phys.}, {\bf 108}:489--526, (1987).

\bibitem{AN84}
M.~Aizenman and C.M. Newman.
\newblock Tree graph inequalities and critical behavior in percolation models.
\newblock {\em J. Stat. Phys.}, {\bf 36}:107--143, (1984).

\bibitem{AKS82}
M.~Ajtai, J.~Koml\'os, and E.~Szemer\'edi.
\newblock Largest random component of a $k$-cube.
\newblock {\em Combinatorica}, {\bf 2}:1--7, (1982).

\bibitem{ABS02}
N.~Alon, I.~Benjamini, and A.~Stacey.
\newblock Percolation on finite graphs and isoperimetric inequalities.
\newblock To appear in {\em Ann.\ Probab.}

\bibitem{AS00}
N.~Alon and J.H. Spencer.
\newblock {\em The Probabilistic Method}.
\newblock Wiley, New York, 2nd edition, (2000).

\bibitem{BK94}
B.~Bollob\'as and Y.~Kohayakawa.
\newblock Percolation in high dimensions.
\newblock {\em Europ. J. Combinatorics}, {\bf 15}:113--125, (1994).

\bibitem{BKL92}
B.~Bollob\'as, Y.~Kohayakawa, and T.~{\L}uczak.
\newblock The evolution of random subgraphs of the cube.
\newblock {\em Random Struct.\ Alg.}, {\bf 3}:55--90, (1992).

\bibitem{BCHSS04a}
C.~Borgs, J.T. Chayes, R.~van~der Hofstad, G.~Slade, and J.~Spencer.
\newblock Random subgraphs of finite graphs: {I}. {The} scaling window under
  the triangle condition.
\newblock Preprint, (2003).

\bibitem{BCHSS04b}
C.~Borgs, J.T. Chayes, R.~van~der Hofstad, G.~Slade, and J.~Spencer.
\newblock Random subgraphs of finite graphs: {II}. {The} lace expansion and the
  triangle condition.
\newblock To appear in {\em Ann. Probab}.

\bibitem{BCHSS04c}
C.~Borgs, J.T. Chayes, R.~van~der Hofstad, G.~Slade, and J.~Spencer.
\newblock Random subgraphs of finite graphs: {III}. {The} phase transition for
  the $n$-cube.
\newblock Preprint, (2003).

\bibitem{CD83}
J.T. Cox and R.~Durrett.
\newblock Oriented percolation in dimensions $d \geq 4$: bounds and asymptotic
  formulas.
\newblock {\em Math. Proc. Cambridge Philos. Soc.}, {\bf 93}:151--162, (1983).

\bibitem{GR78}
D.S. Gaunt and H.~Ruskin.
\newblock Bond percolation processes in $d$ dimensions.
\newblock {\em J. Phys. A: Math. Gen.}, {\bf 11}:1369--1380, (1978).

\bibitem{Gord91}
D.M. Gordon.
\newblock Percolation in high dimensions.
\newblock {\em J. London Math. Soc. (2)}, {\bf 44}:373--384, (1991).

\bibitem{Grim99}
G.~Grimmett.
\newblock {\em Percolation}.
\newblock Springer, Berlin, 2nd edition, (1999).

\bibitem{HS90a}
T.~Hara and G.~Slade.
\newblock Mean-field critical behaviour for percolation in high dimensions.
\newblock {\em Commun. Math. Phys.}, {\bf 128}:333--391, (1990).

\bibitem{HS93up}
T.~Hara and G.~Slade.
\newblock Unpublished appendix to \cite{HS95}. {Available} as paper~93-288 at
  {\tt http://www.ma.utexas.edu/mp\underline{~}arc}.
\newblock (1993).

\bibitem{HS95}
T.~Hara and G.~Slade.
\newblock The self-avoiding-walk and percolation critical points in high
  dimensions.
\newblock {\em Combin. Probab. Comput.}, {\bf 4}:197--215, (1995).

\bibitem{HS04b}
R.~van~der Hofstad and G.~Slade.
\newblock Asymptotic expansions in $n^{-1}$ for percolation critical values on
  the $n$-cube and ${\mathbb Z}^n$.
\newblock Preprint, (2003).

\bibitem{Kest90}
H.~Kesten.
\newblock Asymptotics in high dimensions for percolation.
\newblock In G.R. Grimmett and D.J.A. Welsh, editors, {\em Disorder in Physical
  Systems}. Clarendon Press, Oxford, (1990).

\bibitem{Mens86}
M.V. Menshikov.
\newblock Coincidence of critical points in percolation problems.
\newblock {\em Soviet Mathematics, Doklady}, {\bf 33}:856--859, (1986).

\end{thebibliography}
